\documentclass[12pt]{amsart}

\usepackage[mathscr]{eucal}
\usepackage{amsmath,amssymb,amscd,amsthm}
\usepackage{color}
\usepackage{epsfig}
\newtheorem{Theorem}{Theorem}
\newtheorem{Lemma}{Lemma}

\newtheorem{Definition}{Definition}
\newtheorem{Corollary}{Corollary}

\newtheorem{Proposition}{Proposition}
\theoremstyle{definition}
\newtheorem{Remark}{Remark}

\theoremstyle{plane}

\def \beq{ \begin{equation} }
\def \eeq{\end{equation}}

\title{An intrinsic approach in the curved $n$-body 
    problem: the negative curvature case}
    
\begin{document}

\maketitle

\markboth{F.\ Diacu, E.\ P\'erez-Chavela, and J.G.\ Reyes Victoria}{An intrinsic approach in the curved $n$-body problem of negative curvature}
    
\vspace{-0.5cm}

\author{
\begin{center}
{\rm Florin Diacu \\
Pacific Institute for the Mathematical Sciences\\
and\\
         Departament of Mathematics and Statistics  \\
         University of Victoria \\
         P.O. Box 3060 STN CSC \\
         Victoria, BC, Canada, VSW 3R4\\
         {\tt diacu@uvic.ca}\\
        \medskip
          Ernesto P\'erez-Chavela \\
         Departamento de Matem\'aticas \\
         UAM-Iztapalapa \\
         M\'exico, D.F. MEXICO \\
         {\tt epc@xanum.uam.mx}\\
        \medskip
          J. Guadalupe Reyes Victoria\\
         Departamento de Matem\'aticas \\
         UAM-Iztapalapa \\
         M\'exico, D.F. MEXICO \\
         {\tt revg@xanum.uam.mx}}
\end{center}
}
 
 \bigskip        
\begin{center}
\today
\end{center}



\begin{abstract}
We consider the motion of $n$ point particles of positive masses that interact gravitationally on the 2-dimensional hyperbolic sphere, which has negative constant  Gaussian curvature. Using the stereographic projection, we derive the equations of motion of this curved $n$-body problem in the Poincar\'e disk, where we study the elliptic relative equilibria. Then we obtain the equations of motion in the Poincar\'e upper half plane in order to analyze the hyperbolic and parabolic relative equilibria. Using techniques of Riemannian geometry, we characterize each of the above classes of periodic orbits. For $n=2$ and $n=3$ we recover some previously known results and find new qualitative results about relative equilibria that were not apparent in an extrinsic setting.
\end{abstract}

\newpage

\section{Introduction}
\label{sec:intro}

We consider the negative curvature case of the curved $n$-body problem, i.e.\ the motion of $n$ point masses in spaces of constant Gaussian curvature under the influence of a gravitational potential that naturally extends Newton's law. This article is the completion of an earlier study done for positive curvature by Ernesto P\'erez-Chavela and J.\ Guadalupe Reyes Victoria, \cite{Perez}. The novelty of the approach used in these twin papers is the introduction of intrinsic coordinates, which until recently seemed to be out of reach due to the elaborated computations they involve. But after Florin Diacu, Ernesto P\'erez-Chavela, and Manuele Santoprete obtained the general equations of motion for any number $n$ of bodies in terms of extrinsic coordinates, \cite{Diac}, a glimpse of hope appeared for an intrinsic approach, which has been now achieved for the 2-dimensional case in \cite{Perez} and the present paper. An intrinsic study of the 3-dimensional case is yet to be done.

Whereas \cite{Perez} had set a basic strategy, which we followed here too, in the present study we met with more technical difficulties. In \cite{Perez}, the idea to use stereographic projection and analyze the relative equilibria in the spherical plane had been successful, but the similar approach proved not strong enough for negative curvature. Apart from deriving the equations of motion in the Poincar\'e disk to study elliptic relative equilibria, we also had to use the Poincar\'e upper-half-plane model of hyperbolic geometry in order to overcome the hurdles encountered in the disk model when analyzing hyperbolic and parabolic relative equilibria. This outcome seems to confirm an old saying, which claims that ``in celestial mechanics, there is no ideal system of coordinates.'' The results are, of course, model independent, but it appears that each of them is easier to prove and understand with the help of a particular model.

\subsection{Motivation and history}

The motivation behind the study of the curved $n$-body problem runs deep. 
In the early 1820s, Carl Frie\-drich Gauss allegedly tried to determine the nature of the physical space within the framework of classical mechanics, i.e.\ under the assumption that space and time exist apriori and are independent of each other, \cite{Mill}, \cite{Goe}. He measured the angles of a triangle formed by three mountain peaks, hoping to learn whether space is hyperbolic or elliptic. But the results of his measurements did not deviate form the Euclidean space beyond the unavoidable measurement errors, so his experiments were  inconclusive. Since we cannot reach distant stars to measure the angles of large triangles, Gauss's method is of no practical use for astronomic distances either.

But celestial mechanics can help us find a new approach towards establishing the geometric nature of the physical space. If we extend Newton's $n$-body problem beyond the Euclidian case and also prove the existence of solutions that are specific to each of the negative, zero, and positive constant Gaussian curvature spaces, then the problem of understanding the geometry of the universe reduces to finding, in nature, some of the orbits proved mathematically to exist.  

Therefore obtaining the natural extension of the Newtonian $n$-body problem to spaces of non-zero constant Gaussian curvature, and studying the system of differential equations thus derived, appears to be a worthy endeavour towards comprehending the geometry of the gravitational universe. Additionally, an investigation of this system when the curvature tends to zero may help us better understand the dynamics of the classical case, viewed as a particular problem within a general mathematical framework.

The attempts to find a suitable extension of Newton's gravitational law to spaces of constant curvature started in the 1830s with the work of J\'anos Bolyai,
\cite{Bol}, and Nikolai Lobachevsky, \cite{Lob}, who considered a 2-body problem in hyperbolic space. Unfortunately, their purely geometric approach led to no immediate results. An important step forward was made in 1870 by Ernest Schering, who expressed the ideas of Bolyai and Lobachevsky in analytic terms, \cite{Sche}. This achievement led to the formulation of the equations of motion of the 2-body problem, both for negative and positive curvature. Since then, many researchers brought contributions to the curved 2-body problem, including Wilhelm Killing, \cite{Kil1}, \cite{Kil2}, \cite{Kil3}, and Heinrich Liebmann, \cite{Lie1}, \cite{Lie2}, \cite{Lie3}. The latter showed that the orbits of the curved Kepler problem (i.e.\ the motion of a body around a fixed center) are conics in the hyperbolic plane and proved an analogue of Bertrand's theorem, \cite{Ber}, which states that for the Kepler problem all bounded orbits are closed. But the most convincing argument that the  extension of Newton's law due to Bolyai and Lobachevsky is natural rests with the fact that the potential of the Kepler problem is a harmonic function in the 3-dimensional space, i.e.\ a solution of the Laplace-Beltrami equation, in analogy with the classical potential, which is a solution of Laplace's equation, i.e.\ it is a harmonic function in the classical sense, \cite{Kozlov}. 

Recently, Florin Diacu, Ernesto P\'erez-Chavela, and Manuele Santoprete proposed a new setting for the problem, which allowed an easy derivation of the equations of motion for any $n\ge 2$ in terms of extrinsic coordinates, \cite{Diac}. The combination of two main ideas helped them achieve this goal: the use of Weierstrass's hyperboloid model of hyperbolic geometry (also called the hyperbolic sphere) and the application of the variational approach of constrained Lagrangian dynamics, \cite{Gel}. They also succeeded to solve Saari's conjecture in the collinear case for both positive and negative curvature, having settled the Euclidean case earlier, \cite{Diacu4}, \cite{Diacu5}. Remarkable is also their discovery of the fact that, in the curved 3-body problem, all Lagrangian orbits (i.e.\ rotating equilateral triangles) must have equal masses. Since Lagrangian orbits of non-equal masses exist in our solar system, such as the rotating equilateral triangles formed by the Sun, Jupiter, and the Trojan asteroids, it means that, at distances of the order of $10^1$ AU, space is Euclidean.  

Since then, several papers analyzed the equations of motion. Florin Diacu studied the singularities of the equations and of the solutions, \cite{Diacu2}, Florin Diacu and Ernesto P\'erez-Chavela provided a complete classification of the homographic orbits in the curved 3-body problem, \cite{DiacuPerez}, Florin Diacu considered the polygonal homographic orbits for any finite number of bodies, \cite{Diacu1}, and Pieter Tibboel solved an open problem stated in \cite{Diacu1}. Along the same lines, Regina Mart\'inez and Carles Sim\'o studied the stability of the Lagrangian relative equilibria and homographic orbits for the unit sphere, \cite{Simo}. All these papers treat the 2-dimensional case. The only study, so far, of the 3-dimensional curved $n$-body problem is \cite{Diacu3}, in which Florin Diacu analyzed relative equilibria, a paper in which more details on the history of the problem, as well as an extensive bibliography, can be found.

\subsection{Our results}

As mentioned earlier, this paper is a natural continuation and completion of \cite{Perez}, which was a study of the curved $n$-body problem with the help of intrinsic coordinates
in 2-dimensional spaces of positive curvature. Here we consider the negative curvature case. 

In section 2, we start with the extrinsic description of the motion of $n$ point particles of positive masses on the hyperbolic sphere $\mathbb{L}_R^2$, given by 
the upper sheet of the hyperboloid of 2 sheets,
\begin{equation}\label{upper-sheet}
x^2+y^2-z^2= -R^2,\ \ z>0,
\end{equation}
embedded in the Minkowski space $\mathbb{R}^3_1$. Using stereographic projection through the north pole of the hyperbolic sphere $\mathbb{L}^2_R$, we move the problem to the Poincar\'e disk $\mathbb{D}^2_R$, where we obtain the equations of motion in complex coordinates and show that the original equations of motion are equivalent to the equations obtained in $\mathbb{D}^2_R$. Then we provide the first integrals of the system in $\mathbb{D}^2_R$.

In section \ref{sec:rela-equilibria}, we consider a suitable Killing vector field
in $\mathbb{D}^2_{R}$ and its associated one-dimensional additive subgroup of the Lie algebra $\mathfrak{su}(1,1)$. When we project this Lie algebra via the exponential map onto the Lie group $SU(1,1)$, we obtain a one-dimensional subgroup of proper isometries which generate the elliptic  relative equilibria of the curved $n$-body problem in $\mathbb{D}^2_{R}$. By analyzing the corresponding one-dimensional subgroup of M\"obius transformations, we obtain the equations that characterize all the elliptic
relative equilibria. Then we study these elliptic relative equilibria for $n=2$ and $n=3$.
In the former case we give a complete description of the elliptic relative equilibria. In the latter case we study both the Eulerian and the Lagrangian orbits and describe their qualitative behavior, thus recovering in $\mathbb D_R^2$ some results proved in \cite{DiacuPerez} for $\mathbb L_R^2$. Additionally, we show that other elliptic relative equilibria don't exists in the curved 3-body problem in the case of negative curvature.

Although we also obtain in section 3 the one-dimensional subgroups of proper isometries that generate the hyperbolic and parabolic relative equilibria, their complicated expressions make these orbits hard to understand. Therefore, in section \ref{sec:Klein_model}, we move the problem from the Poincar\'e disk $\mathbb D_R^2$ to the Poincar\'e upper half plane $\mathbb{H}^2_R$, in which we obtain the intrinsic equations of motion. Again, we consider the suitable Killing vector fields in $\mathbb{H}^2_R$  and their corresponding one-parameter additive subgroups in its Lie Algebra $\mathfrak{sl}(2,\mathbb{R})$, which we project via the exponential map onto the Lie group of proper isometries $SL(2,\mathbb{R})$. The subgroups corresponding to elliptic relative equilibria become complicated in this context, but we already studied them in the Poincar\'e disk. Fortunately, the subgroups corresponding to hyperbolic and parabolic relative equilibria are simple enough to allow a complete analysis. So using the associated one-parameter subgroups of M\"obius transformations, we obtain the conditions for the existence of hyperbolic and parabolic relative equilibria. Then we proceed as in section 3 with an analysis of the cases $n=2$ and $n=3$ for the former orbits, whose qualitative behavior we describe if they are of Eulerian type. Finally we confirm, within our more general context, a result obtained in  \cite{Diac}, which proves that parabolic relative equilibria do not exist.


\section{Equations of motion}\label{sec:equamotion}

Consider a connected and simply connected 2-dimensional surface of
constant negative Gaussian curvature $\displaystyle \kappa=-1/{R^2}$.
It is well known (see, e.g., \cite{DoCarmo}) that this surface
can be locally represented by \eqref{upper-sheet}, the upper sheet of the 2-dimensional
hyperbolic sphere of radius $iR$, denoted by $\mathbb{L}^2_R$, embedded in the Minkowski space $\mathbb{R}^3_1$, which is endowed with 
the Lorentz inner product defined by
\begin{equation}\label{Lorentz}
Q_1 \odot Q_2:= {\mathfrak x_1\mathfrak x_2+\mathfrak y_1 \mathfrak y_2-\mathfrak z_1 \mathfrak z_2}
\end{equation}
for any points $Q_1=(\mathfrak x_1,\mathfrak y_1,\mathfrak z_1)$ and $Q_2=(\mathfrak x_2,\mathfrak y_2,\mathfrak z_2)$ of $\mathbb{R}^3_1$.

Let $\Pi$ be the stereographic projection of $\mathbb{L}^2_R$
through the north pole, $(0,0, -R)$, located at the vertex of the lower sheet of
the hyperbolic sphere. The image of $\Pi$ is the planar disk of radius $R$, denoted by $\mathbb{D}^2_R$, having the center at the origin of coordinate system, so we can write this function as 
\begin{eqnarray}\label{stereo}
\Pi \colon \mathbb{L}^2_R \to \mathbb{D}^2_R, \ \ \Pi(Q) = q. 
\end{eqnarray}
It is easy to see that $\Pi(\mathbb{L}^2_R)$ is the entire open disk
 $\mathbb{D}^2_R$, and that if we take $Q=(\mathfrak x,\mathfrak y,\mathfrak z)$, then $q=(u,v)$, where
\begin{equation} \label{ste-eq}
   u= \frac{R\mathfrak x}{R+\mathfrak z} \quad  \quad {\rm and}  \quad  \quad v= \frac{R\mathfrak y}{R+\mathfrak z}.
\end{equation}
The inverse of the stereographic projection is given by the equations
\begin{equation}\label{ste-inv}
   \mathfrak x= \frac{2R^2u}{R^2-u^2-v^2}, \quad  \mathfrak y= \frac{2R^2v}{R^2-u^2-v^2}, \quad \mathfrak z
= \frac{R(R^2+u^2+v^2)}{R^2-u^2-v^2}.
\end{equation}
The metric (distance) of the sphere  $\mathbb{L}_R^2$ is transformed
into the metric
\begin{equation}\label{met}
   ds^2= \frac{4R^4}{(R^2-u^2-v^2)^2} (du^2 + dv^2).
\end{equation}
In the Poincar\'e disk $\mathbb{D}^2_R $, the geodesics are the diameters of the circle
of radius $R$ and the arcs of the circles orthogonal to it. In terms of the above metric, all geodesics have infinite length, so the circle of radius $R$ represents the points at infinity.

Since the metric is conformal, with factor of conformity 
$$
\displaystyle\lambda(u,v)=\frac{4R^4}{(R^2-u^2-v^2)^2},
$$ 
the Christoffel symbols associated to the metric are given by
\begin{equation}\label{chrishyp}
 -\Gamma_{22}^1 = \Gamma_{11}^2=\Gamma_{12}^1=\frac{1}{2 \lambda(u,v)}\frac{\partial \lambda}{\partial u}=
 \frac{2 u}{(R^2-u^2-v^2)^2}, 
 \end{equation}
 \begin{equation}
 \Gamma_{22}^2 = -\Gamma_{11}^1=\Gamma_{12}^2=\frac{1}{2 \lambda(u,v)}\frac{\partial \lambda}{\partial v}=
 \frac{2 v}{(R^2-u^2-v^2)^2}, 
 \label{christoffelhyp}
\end{equation} 
so the geodesics can also be obtained by solving the system of second order differential equations (see \cite{Dub} for more details):
\begin{eqnarray*}
\ddot{u} + \Gamma_{11}^1 \dot{u}^2 +2 \Gamma_{12}^1 \dot{u} \dot{v}
+\Gamma_{22}^1 \dot{v}^2 & = & 0, \\[0.2cm]
\ddot{v} + \Gamma_{11}^2 \dot{u}^2 +2 \Gamma_{12}^2 \dot{u} \dot{v}
+\Gamma_{22}^2 \dot{v}^2 & = &0,
\end{eqnarray*}
which is equivalent to 
\begin{equation}
\label{geohyp}
\begin{cases}
\ddot{u} = - \frac{2}{R^2-u^2-v^2}\left( u\dot{u}^2 + 2v\dot{u}\dot{v} -
  u\dot{v}^2 \right)\cr
    \ddot{v} = - \frac{2}{R^2-u^2-v^2}\left( v\dot{v}^2 + 2u\dot{u}\dot{v} -
  v\dot{u}^2 \right).\cr 
  \end{cases}
\end{equation}
From now on, we will think of $ \mathbb{D}^2_{R}$ as being endowed with the above metric. The geodesic distance between
two points  $q_k, q_j \in \mathbb{D}^2_{R}$ is
$$ d_{kj}= d(q_k,q_j) = d(Q_k, Q_j )= R \cosh^{-1} \left(-{\frac{Q_k \odot Q_j}{R^2}}\right).$$

In analogy with \cite{Perez}, where we considered the cotangent potential for $\kappa>0$, we will use here the hyperbolic cotangent potential for $\kappa<0$. 
With the help of the stereographic projection, we see that
$$ 
Q_k \odot Q_j = \frac{4R^4 \,  q_k \cdot q_j -
R^2(||q_k||^2 + R^2) (||q_j||^2 + R^2) }{(R^2-||q_k||^2  )(R^2-||q_j||^2
)},
$$ 
where $\cdot$ is the standard inner product in the
plane $\mathbb{R}^2$. Then
\begin{equation}\label{cot-r}
  \coth_R{\left(\frac{d_{kj}}{R}\right)}= -\frac{4R^2 \, q_k \cdot q_j - (||q_k||^2 + R^2)( (||q_j||^2 + R^2)}{W},
 \end{equation}
 where $$
W=  {\sqrt{W_1 - W_2}},$$
$$W_1 =  \left[ 4R^2 \,  q_k \cdot q_j - (||q_k||^2 + R^2)( (||q_j||^2 + R^2)\right]^2,$$
$$W_2 = (R^2-||q_k||^2 )^2((R^2-||q_j||^2 )^2.$$

We now introduce complex variables in $\mathbb{D}^2_{R}$, 
$$
z=u+iv,\ \ \bar{z}=u-iv,
$$ 
in order to simplify the computations. The inverse of this transformation is given by
\begin{equation}\label{invz}
u=\frac{z+\bar{z}}{2},\ \ v=\frac{z-\bar{z}}{2i}.
\end{equation}
Then in terms of $z$ and $\bar z$, the formulas, \eqref{ste-inv}, for the inverse of the stereographic projection are given by
\begin{equation}\label{inv-z-stereo}
\mathfrak x=\frac{R^2(z+\bar z)}{R^2-|z|^2},\ \ \mathfrak y=\frac{iR^2(\bar z-z)}{R^2-|z|^2},\ \ \mathfrak z=\frac{R(|z|^2+R^2)}{R^2-|z|^2},
\end{equation}
where $|.|$ denotes the absolute value of a complex number.

The distance (\ref{met}) and equation (\ref{cot-r})
take the form
\begin{equation}\label{met-c}
   ds^2= \frac{4R^4}{(R^2-|z|^2)^2}\, dz d\bar{z}
\end{equation}
and
\begin{equation}\label{eq:coth}
 \coth_R \left(\frac{d_{kj}}{R}\right)=-\frac{2(z_k \bar{z}_j+z_j \bar{z}_k)R^2-
(|z_k|^2+R^2)(|z_j|^2+R^2)}{\sqrt{\Theta_{2,(k,j)}(z, \bar{z})}}, \\
\end{equation}
where
\begin{equation}\label{eq:discoth}
\Theta_{2,(k,j)}(z, \bar{z})=[2(z_k \bar{z}_j+z_j \bar{z}_k)R^2
-(|z_k|^2+R^2)(|z_j|^2+R^2)]^2
\end{equation}
$$
-(R^2-|z_k|^2)^2(R^2-|z_j|^2)^2
$$
and $d_{kj}=d(z_k,z_j)$ denotes the geodesic distance in the metric
(\ref{met-c}) between the points $z_k$ and $z_j$ in  $\mathbb{D}^2_{R}$. From now on we will think of  $\mathbb{D}^2_{R}$ as the hyperbolic Poincar\'e disk endowed with this new form of the metric.

Notice that, in these complex coordinates, system \eqref{geohyp}, which yields the geodesics, takes the form
\begin{equation} \label{geodesics}
  \ddot{z}+\frac{ 2 \bar{z} \dot{z}^2}{ R^2- |z|^2}=0.
\end{equation}

\subsection{The intrinsic approach}

Let  $z=(z_1,z_2,\cdots,z_n) \in (\mathbb{D}^2_R)^n$ be the configuration of
$n$ point particles with masses $m_1,m_2,\cdots,m_n>0$ in $\mathbb{D}^2_{R}$.  We assume that the particles are moving
under the action of the Lagrangian
\begin{equation}\label{eq:lagrangiangeral}
 L_R \, (z,\dot{z},\bar{z},\dot{\bar{z}})= T_R(z,\dot{z},\bar{z},\dot{\bar{z}}) + U_R (z, \bar{z}),
\end{equation}
 where
 \begin{equation}\label{eq:kinetic-energy}
T_R \, (z,\dot{z},\bar{z},\dot{\bar{z}}) = \frac{1}{2}\sum_{1 \leq k \leq n} m_k \, \lambda (z_k, \bar{z}_k) \, |\dot{z}_k|^2
\end{equation}
is the kinetic energy and
\begin{equation}\label{eq:potesf}
U_R (z, \bar{z}) = \frac{1}{R} \sum_{1 \leq k < j \leq n}^n m_k m_j \coth_R
\left(\frac{d_{kj}}{R}\right)
\end{equation}
$$
=- \frac{1}{R} \sum_{1 \leq k < j \leq n}^n m_k m_j \frac{2(z_k \bar{z}_j+z_j \bar{z}_k)R^2-
(|z_k|^2+R^2)(|z_j|^2+R^2)}{\sqrt{\Theta_{2,(k,j)}(z, \bar{z})}}
$$
is the force function of the particle system ($-U_R$ is the potential) defined in the set $\mathbb{D}^{2n}_{R}
\setminus \Delta$, where $\Delta$ is the set of zeroes of $\Theta_{2,(k,j)}(z, \bar{z})$, given by (\ref{eq:discoth}), 
and
\begin{equation}\label{eq:conforesf}
\lambda (z_k, \bar{z}_k)= \frac{4R^4}{(R^2-|z_k|^2)^2}
\end{equation}
is the conformal function of the Riemannian metric.

The following result gives  the equations of motion for the problem. It's proof is perfectly similar with the homologue result for curvature $\kappa>0$, obtained in \cite{Perez}, so we omit it here.

\begin{Lemma}\label{lema:gral-mechanicalhyp}  Let
\[ \displaystyle L(z,\dot{z})= \frac{1}{2}\sum_{k=1}^n m_k \lambda (z_k, \bar{z}_k) \, |\dot{z}_k|^2 + U_R(z, \bar{z})  \]
be the Lagrangian defined in \eqref{eq:lagrangiangeral} for the
given problem. Then the  solutions of the corresponding
Euler-Lagrange equations satisfy the system of second order
differential equations
\begin{equation} \label{eq:motiondisk}
  m_k  \ddot{z}_k = -\frac{ 2 m_k \bar{z_k} \dot{z}_k^2}{ R^2- |z_k|^2  } + \frac{2}{\lambda(z_k, \bar{z}_k)} \, \frac{\partial U_R}{\partial \bar{z}_k},\ k=1,\dots, n,
\end{equation}
where
\begin{equation}\label{eq:gradmet2}
 \frac{\partial U_R}{\partial \bar{z}_k}= \sum_{\substack{j =1\\ j\ne k}}^n \frac{2 m_k m_j R P_{2,(k,j)} (z, \bar{z})}{(\Theta_{2,(k,j)}(z, \bar{z}))^{3/2}},
\end{equation}
$$ 
P_{2,(k,j)}(z, \bar{z})=(R^2-|z_k|^2)(R^2-|z_j|^2)^2(z_j - z_k)(R^2 - z_k \bar{z}_j), 
$$
$$
\Theta_{2,(k,j)}(z, \bar{z})=[2(z_k \bar{z}_j+z_j \bar{z}_k)R^2
-(|z_k|^2+R^2)(|z_j|^2+R^2)]^2
$$
$$
-(R^2-|z_k|^2)^2(R^2-|z_j|^2)^2,\ k,j\in\{1,\dots,n\}, k\ne j.
$$
\end{Lemma}

Notice that the first term on the right hand side of \eqref{eq:motiondisk} depends on
the kinetic energy alone, whereas the second term depends on the potential.
Therefore, from Lemma \ref{lema:gral-mechanicalhyp}, we can also draw the following conclusion.

\begin{Corollary}
\label{sec:free} If in $\mathbb{D}^2_{R}$ the potential is constant in the entire space, then the particles move freely along geodesics.
\end{Corollary}

\begin{proof}
Since the potential is constant, equations (\ref{eq:motiondisk}) take the form
\begin{equation} \label{eq:motiondiskgeo}
  m_k  \ddot{z}_k + \frac{ 2 m_k \bar{z}_k \dot{z}_k^2}{ R^2- |z_k|^2  } = 0,\ k=1,\dots, n,
\end{equation}
and, if we simplify by $m_k>0$, we obtain that the coordinates of each body satisfy equation \eqref{geodesics},  so each body moves along a geodesic. 
\end{proof}

\subsection{Equivalence of the models}

In \cite{Diac}, where $\widetilde\nabla=(\partial_{\mathfrak x},\partial_{\mathfrak y},-\partial_{\mathfrak z})$, the equations of motion for the $n$-body problem in
the hyperbolic space $\mathbb{L}^2_R$
are
\begin{equation}
m_j \ddot Q_j=\widetilde\nabla_{Q_j} \mathcal V_\kappa(Q)
-m_j \kappa(\dot Q_j\odot\dot Q_j)Q_j, 
 \ \ j=1,\dots, n,\label{eqmotiondisk}
\end{equation}
where $Q=(Q_1,\dots, Q_n)$ is the configuration of the system, $\mathcal V_\kappa$ is the force function, 
\begin{equation}
{\mathcal V}_\kappa(Q)=\sum_{1\le l<j\le n}{m_lm_j
|\kappa|^{1/2}{\kappa{Q}_l\odot Q_j}\over
[({\kappa{Q}_i\odot{Q}_j})^2-(\kappa{Q}_l
\odot{Q}_l)(\kappa{Q}_j\odot{Q}_j)]^{1/2}},
\end{equation}
$Q_j\odot Q_j=\kappa^{-1}=R^2,\ j=1,\dots,n$, are the constraints, which maintain the particles on the hyperbolic sphere, and $\odot$ denotes the Lorentz inner product of the Minkowski space, defined in \eqref{Lorentz}.
The $Q_l$-gradient, $l=1,\dots,n$, of the force function is
\begin{equation}
\widetilde\nabla_{Q_l}\mathcal V_\kappa(Q)=\sum_{j=1,j\ne l}^n{{m_lm_j}|\kappa|^{3/2}
\left[Q_j  - (\kappa Q_l\odot Q_j)Q_l \right]
\over
\left[\left(\kappa Q_l\odot Q_j\right)^2-1\right]^{3/2}}.
\label{grad}
\end{equation}
As shown in \cite{Diac}, system \eqref{eqmotiondisk} has 4 integrals of motion. One is the integral of energy,
\begin{equation}\label{Leng}
\frac{1}{2}\sum_{j=1}^nm_j(\dot{Q}_j\odot\dot{Q}_j)(\kappa Q_j\odot Q_j)-{\mathcal V_\kappa}(Q)=h,
\end{equation}
where $h$ is the energy constant, and the other 3 are the integrals of the total angular momentum,
\begin{equation}\label{Lang1}
\sum_{j=1}^nm_j({\mathfrak y}_j\dot{\mathfrak z}_j-{\mathfrak z}_j\dot{\mathfrak y}_j)=c_1,
\end{equation}
\begin{equation}\label{Lang2}
\sum_{j=1}^nm_j({\mathfrak z}_j\dot{\mathfrak x}_j-{\mathfrak x}_j\dot{\mathfrak z}_j)=c_2,
\end{equation}
\begin{equation}\label{Lang3}
\sum_{j=1}^nm_j({\mathfrak y}_j\dot{\mathfrak x}_j-{\mathfrak x}_j\dot{\mathfrak y}_j)=c_3,
\end{equation}
where $c_1, c_2, c_3$ are the constants of the total angular momentum.

Now we can state the equivalence between the equations of motion of
the curved $n$-body problem defined on the hyperbolic sphere $\mathbb{L}_R^2$, of
Gaussian curvature $\kappa=-1/R^2<0$, and the equations we have obtained in the Poincar\'e disk $\mathbb{D}^2_{R}$. The proof of this result is perfectly similar with the proof
of Theorem 2.3 in \cite{Perez}, which concerns the case $\kappa>0$, so we omit it here.

\begin{Theorem}
\label{equiv} The equations of motion of the curved $n$-body problem  on
the Poincar\'e disk $\mathbb{D}^2_{R}$, given by system \eqref{eq:motiondisk}, and the
corresponding equations on the hyperbolic sphere  $\mathbb{L}_R^2$, given by system \eqref{eqmotiondisk}, are equivalent.
\end{Theorem}

\subsection{Integrals of motion}

Since Theorem \ref{equiv} proves the equivalence between the equations of motion of the curved $n$-body problem in $\mathbb L_R^2$ and $\mathbb D_R^2$, system \eqref{eq:motiondisk} inherits the properties of system \eqref{eqmotiondisk}. We are interested here in how the first integrals of motion translate from system  \eqref{eqmotiondisk} to system \eqref{eq:motiondisk}. Using equations \eqref{inv-z-stereo}, which give the inverse of the stereographic projection, we see that the integrals of motion given by \eqref{Leng}, \eqref{Lang1}, \eqref{Lang2}, and \eqref{Lang3} become, respectively,
\begin{equation}\label{Deng}
\frac{1}{2}\sum_{k=1}^n m_k \lambda (z_k, \bar{z}_k) |\dot{z}_k|^2 - U_R(z, \bar{z})=h,
\end{equation}
where $h$ is the same energy constant as in \eqref{Leng}, and
\begin{equation}\label{Dang1}
\sum_{k=1}^n\frac{im_kR^3}{(R^2-|z_k|^2)^2}[R^2(\dot{z}_k-\dot{\bar{z}}_k)+\dot{z}_k\bar{z}_k^2-\dot{\bar{z}}_kz_k^2]=c_1,
\end{equation}
\begin{equation}\label{Dang2}
\sum_{k=1}^n\frac{m_kR^3}{(R^2-|z_k|^2)^2}[R^2(\dot{z}_k+\dot{\bar{z}}_k)-\dot{z}_k\bar{z}_k^2-\dot{\bar{z}}_kz_k^2]=c_2,
\end{equation}
\begin{equation}\label{Dang3}
\sum_{k=1}^n\frac{2im_kR^4}{(R^2-|z_k|^2)^2}(\bar{z}_k\dot{z}_k-z_k\dot{\bar{z}}_k)=c_3.
\end{equation}
A straightforward computation confirms that the left hand sides of equations \eqref{Dang1}, \eqref{Dang2}, and \eqref{Dang3} are real functions. Moreover, $c_1, c_2, c_3\in\mathbb R$ are the same constants as in \eqref{Lang1}, \eqref{Lang2}, and \eqref{Lang3}.

\section{Relative equilibria in $\mathbb{D}^2_{R}$}\label{sec:rela-equilibria}

In this section we start to analyze the dynamics of the particles that interact in  $\mathbb{D}^2_{R}$. Our goal is to give some general characterization of the relative equilibria, which we define below with the help of geometric tools. For this purpose, we will need to understand first what are the singularities of the equations of motion, such that to avoid the singular configurations when dealing with regular equilibria, and to introduce the Principal Axis Theorem, which will guide us towards finding a proper definition of the relative equilibria.

\subsection{The singularities of the equations of motion in  $\mathbb{D}^2_R$ }

The equations of motion \eqref{eq:motiondisk} have singularities in $\mathbb{D}^2_R$  when at least one denominator vanishes, i.e.\ if there exist indices $k,j$, with $1\le k<j\le n$, such that
 \[ \Theta_{2,(k,j)} (z, \bar{z})=0.  \]
A straightforward computation shows that this statement is equivalent to saying that there are $k,j$, with $1\le k<j\le n$, such that
$$
(R^2-z_\kappa\bar{z}_j)(R^2-z_j\bar{z}_\kappa)(\bar{z}_\kappa-\bar{z}_j)(z_j-z_\kappa)=0.
$$
Since $|z_i|,|\bar{z}_i|<R,\ i=1,\dots,n$, it follows that this equation is satisfied only for $z_j=z_k$, i.e.\ when a collision takes place. Therefore the set of singularities of the equations of motion is
$$
\Delta = \bigcup_{1\le k<j\le n} \,  \Delta_{kj},
\ {\rm where}\  \ \Delta_{kj}=\{(z_1,\dots, z_n)\ \! |\ \! z_k = z_j \}.
$$

\subsection{Principal axis theorem in $\mathbb L_R^2$}
\label{subsec:relativesphere}

Let $SO(1,2)$ denote the Lorentz group, defined as the Lie group of all the isometric rotations of determinant 1 in $\mathbb{R}^3_1$ that keep the hyperbolic sphere $\mathbb{L}^2_R$ invariant.  An important result related to this subgroup of orthogonal transformations is the Principal Axis Theorem, \cite{Nomizu}, which states that there are 3 1-parameter subgroups of $SO(1,2)$ whose elements can be represented, in some suitable basis of $\mathbb{R}^3_1$, as
\begin{enumerate}
\item matrices of the form
\begin{equation}
A=P\begin{bmatrix}
\cos\theta & -\sin\theta & 0 \\
\sin\theta & \cos\theta & 0 \\
0 & 0 & 1
\end{bmatrix} P^{-1},
\label{elliptic-rot}
\end{equation}
called {\it elliptic rotations} in $\mathbb{L}^2_R$ around a {\it timelike} axis, the $\mathfrak z$-axis in our case, 
\item matrices of the form
\begin{equation}
B=P\begin{bmatrix}
1 & 0 & 0 \\
0 & \cosh s & \sinh s \\
0 & \sinh s & \cosh s
\end{bmatrix}P^{-1},
\end{equation}
called {\it hyperbolic rotations} in $\mathbb{L}^2_R$, around a {\it spacelike} axis, the $\mathfrak x$-axis in our case, and
\item matrices of the form
\begin{equation}
C=P\begin{bmatrix}
1 & -t & t \\
t & 1-t^2/2 & t^2/2 \\
t & -t^2/2 & 1+t^2/2
\end{bmatrix}P^{-1},
\end{equation}
called {\it parabolic rotations} in $\mathbb{L}^2_R$, around a {\it lightlike} axis, which is given here by $\mathfrak x=0$, $\mathfrak{y=z}$, where $P\in SO(1,2)$ in all cases.
\end{enumerate}

\begin{Remark}  This result implies that any isometry in $\mathbb{L}_R^2$ can be written as a composition of an elliptic rotation around the $\mathfrak z$-axis, a hyperbolic rotation around the $\mathfrak x$-axis, and a parabolic rotation around the axis $\mathfrak y=0, \mathfrak{z=x}$. 
\end{Remark}

\subsection{Relative equilibria  in  $\mathbb{D}^2_{R}$}
\label{subsec:relativesphere}

Let ${\rm Iso}(\mathbb{D}^2_{R})$ be  the group of isometries of $\mathbb{D}^2_{R}$, and assume that $\{G(t)\}$ is a 1-parameter subgroup of ${\rm Iso}(\mathbb{D}^2_{R})$ that leaves $\mathbb{D}^{2n}_{R} \setminus \Delta$ and $\Delta$ invariant. We can now give a general definition for relative equilibria in $\mathbb{D}^2_{R}$. 

\begin{Definition}\label{def:equilibria}
A  {\it relative equilibrium} of the negatively curved $n$-body problem is a solution $z$ of equations \eqref{eq:motiondisk} that is invariant relative to some subgroup $\{G(t) \}$ of ${\rm Iso}(\mathbb{D}^2_{R})$, i.e.\ the function $w$, given by $w(t)= G(t)z(t)$, obtained by the action of some element of $G$, is also a solution of system \eqref{eq:motiondisk}.
\end{Definition}

To understand the implications of this definition and be able to represent the relative
equilibria in $\mathbb D_R^2$ in as precise terms as we did for $\mathbb L_R^2$ in \cite{Diac}, we need to take first a look at the topological group structure of isometric rotations of $\mathbb D_R^2$ (for more details, see, e.g., \cite{Dub}). For this purpose, consider the matrix
$$
   \tilde{I}= \left(\begin{array}{ccc}
    1   &  0     \\
   0    & -1   \\
    \end{array}\right)
$$
and let
$$
{\rm SU}(1,1) = \{ A \in {\rm GL}(2,\mathbb{C}) \, | \, \, \bar{A}^T \tilde{I} A=\tilde{I}, \ \det A=1\} 
$$
be the {\it special orthochronous unitary group}.  Then some algebraic computations show that any matrix $A \in {\rm SU}(1,1)$ has the form
$$
   A= \left(\begin{array}{cc}
    a        &  b     \\
   \bar{b} & \bar{a}   \\
    \end{array}\right),
$$
with $a,b \in \mathbb{C}$ satisfying   $|a|^2 -|b|^2 =1$. This last condition implies that the group ${\rm SU}(1,1)$ is diffeomorphic with the real 3-dimensional unit hyperbolic sphere embedded in $\mathbb{C}^2$. The term {\it orthochronous} means that the transformations $A$ do not change the direction of the time $t$ in the standard interpretation of the Minkowski space. In our case, however, this is just another space coordinate.

The {\it group of proper orthochronous isometries} of $\mathbb{D}^2_R$, i.e.\ the group of transformations that also maintain the geometric orientation, is the quotient group 
$$
\displaystyle {\rm SU}(1,1)/\{\pm I \},
$$
where $I$ is the unit $2\times 2$ matrix. To every class $\displaystyle A \in {\rm SU}(1,1)/\{\pm I \}$, we can associate a M\"obius transformation,  
$$
f_A : \mathbb{D}^2_R \to \mathbb{D}^2_R,\ \ 
f_A (z) = \frac{a z+b}{\bar{b} z + \bar{a}},
$$
for which it is easy to see that $\displaystyle f_{-A} (z)= f_A (z)$.

If $M(2,\mathbb C)$ is the set of $2\times 2$ matrices with complex elements, 
the {\it Lie algebra} of ${\rm SU}(1,1) $ is the 3-dimensional real
linear space
	\[ \mathfrak{su}(1,1) = \{ X \in {\rm M}(2,\mathbb{C}) \, | \, \, \tilde{I} \, \bar{X}^T =-X \,\tilde{I}, \, \, {\rm trace\ \!} X=0  \}  \]
spanned by the Killing vector fields in $\mathbb{D}_R^2$ associated
to the Pauli matrices, 
$$
g_1 = \frac{1}{2} \left( \begin{array}{ccc}
    0 & 1 \\
    1 & 0  \\
    \end{array}\right),  \quad
g_2 = \frac{1}{2}\left(\begin{array}{cc}
    i & 0 \\
    0 & -i  \\
    \end{array}\right), \quad
g_3 = \frac{1}{2}\left(\begin{array}{cc}
    0 & i \\
    -i & 0  \\
    \end{array}\right),
$$
which form a basis of $\mathfrak{su}(1,1)$.

Consider further the exponential map of matrices,
$$
\exp\colon\mathfrak{su}(1,1) \to {\rm SU}(1,1),
$$
applied to  the one-parameter additive subgroups $\{ tg_1 \}$, $ \{ tg_2 \}$, and $\{ tg_3 \}$, which  are straight lines form the geometric point of view. This operation leads us to the following one-parameter subgroups of ${\rm SU}(1,1) $ in $\mathbb D_R^2$:
\begin{enumerate}
\item  the subgroup
\[ G_1(t)=\exp (tg_1)=
\left(\begin{array}{cc}
    \cosh (t/2) & \sinh (t/2) \\
    \sinh (t/2)  & \cosh (t/2)  \\
    \end{array}\right),
\]
which defines the one-parameter family of M\"obius
transformations
\begin{equation} \label{eq:Mobius-Disk-1}
 f_{G_1} (z,t) = \frac{ \cosh (t/2)z   +\sinh (t/2)}{ \sinh (t/2)z   + \cosh (t/2)},
 \end{equation}
 
\item  the subgroup
\[ G_2 (t)=\exp (tg_2)=
\left(\begin{array}{cc}
    e^{it/2} & 0 \\
    0  & e^{-it/2} \\
    \end{array}\right),
\]
which defines the one-parameter family of  M\"obius transformations
\begin{equation} \label{eq:Mobius-Disk-2}
f_{ G_2 } (z,t) = e^{it}z,
\end{equation}

\item  the subgroup
\[ G_3(t)=\exp (tg_3)=
\left(\begin{array}{cc}
    \cosh (t/2) & i\sinh (t/2) \\
    -i\sinh (t/2)  & \cosh (t/2)  \\
    \end{array}\right),
\]
which defines the one-parameter family of M\"obius transformations
\begin{equation} \label{eq:Mobius-Disk-3}
 f_{G_3} (z,t) = \frac{ \cosh (t/2)z   + i \sinh(t/2)}{ - i\sinh (t/2)z   + \cosh (t/2)}.
 \end{equation}
\end{enumerate}

\subsection{Elliptic relative equilibria}
\label{subsec:elliptic-relative eq}

When projected into the Poincar\'e disk, a compositions of the elliptic, hyperbolic, and parabolic rotations of $\mathbb L_R^2$  corresponds to some composition of the M\"obius transformations  (\ref{eq:Mobius-Disk-1}), (\ref{eq:Mobius-Disk-2}), and (\ref{eq:Mobius-Disk-3}).  We therefore need to analyze each of these one-parameter subgroups of transformations in order to find their associated relative equilibria. We will first achieve this goal for the second group, $G_2$, given by the M\"obius transformations $f_{ G_2 } (z) = e^{it}z$, which correspond to the elliptic rotations in $\mathbb L_R^2$. The other 2 groups will be treated in section 4, since their study becomes very tedious in $\mathbb D_R^2$ because of complicated computations.

Notice that the one-parameter rotation subgroup of $SO(1,2)$ introduced
by the matrix
\begin{equation}\label{rot-mat}
   \mathcal A(t) = \left(\begin{array}{ccc}
    \cos{t}        &  -\sin{t}   &  0 \\
   \sin{t} & \cos{t} & 0   \\
     0      &     0   & 1\\
    \end{array}\right),
\end{equation}
i.e.\ the rotation \eqref{elliptic-rot} around the $\mathfrak z$ axis of $\mathbb R_1^3$ that leaves invariant any hyperbolic sphere $\mathbb{L}^2_R$ centered at the origin,
is the isometric flow for the basic {\it Killing vector field}
\begin{equation}\label{eq:killingvector}
L_Z(\mathfrak{x,y,z}) = (\mathfrak{-y,x},0).
\end{equation}
This observation leads us to the following result.

\begin{Proposition}\label{prop:isomorph}
Let $H\colon{\rm SU}(1,1)/\{\pm I \} \to SO(1,2)$ be an isomorphism between the groups of proper orthochronous isometries of the Poincar\'e disk $\mathbb D_R^2$ 
and the Lorentz group of the hyperbolic sphere $\mathbb L_R^2$. Then 
$$
H(G_2(t)) = \mathcal A(t).
$$
\end{Proposition}
\begin{proof}
The stereographic projection,
$$\Pi (\mathfrak{x,y,z}) = \left( \frac{R\mathfrak x}{R+\mathfrak z}, \frac{R\mathfrak y}{R+\mathfrak z} \right),$$
shows that since the rotation tangent vector at $(\mathfrak{x,y,z})$ in
$\mathbb{L}^2_R$ is $(\mathfrak{-y,x},0)$, after the projection we have
\begin{equation}\label{unique}
D\Pi[(L_Z) (\mathfrak{x,y,z})]^T = \left(\begin{array}{ccc}
    \frac{R}{R+\mathfrak z}        &  0  &   \frac{R\mathfrak x}{(R+\mathfrak z)^2} \\
0 &  \frac{R}{R+\mathfrak z} &   \frac{R\mathfrak y}{(R+\mathfrak z)^2}  \\
\end{array}\right)\left(\begin{array}{c}
                                          -\mathfrak y \\ \mathfrak x \\ 0
                                          \end{array}\right) = \left(\begin{array}{c}
                                           -\frac{R\mathfrak y}{R+\mathfrak z} \\ \ \ \  \! \frac{R\mathfrak x}{R+\mathfrak z} \end{array}\right),
 \end{equation}
where $D\Pi$ is the Jacobian matrix and the upper $T$ denotes the transposed of the vector. In complex notation this relationship corresponds to 
$$
\displaystyle -v + iu = i(u + iv) = iz,$$
which leads us to the differential equation
\begin{equation}\label{complex-eq}
\dot{z} = iz.
\end{equation}
Its flow is given by $f_t(z) = e^{it}z$, associated to the one-parameter subgroup of M\"obius transformations $f_{G_2}$. This remark completes the proof. 
\end{proof}

In order to obtain from Proposition \ref{prop:isomorph} some information regarding the relative equilibria of type $G_2$ in $\mathbb{D}^2_{R}$, we consider functions of the form
\begin{equation}\label{expo}
 w_k (t) = e^{it} z_k (t),
 \end{equation}
where $z = (z_1, \dots, z_n)$  is a solution of equation (\ref{eq:motiondisk}), and look for conditions that the function
$
w=(w_1,\dots, w_n)
$
is also a solution of system (\ref{eq:motiondisk}). Straightforward  computations show that
\begin{equation}\label{acting-change}
\begin{cases}
\dot{w}_k=  (i z_k + \dot{z}_k ) e^{it} \cr
\ddot{w}_k=(\ddot{z}_k  +2i \dot{z}_k -z_k ) e^{it} \cr
\frac{d \bar{z}_k}{d \bar{w}_k}= e^{it},\ \ k=1,\dots,n. \cr
\end{cases}
\end{equation}
Using these facts together with the conditions that $w$ is a solution of equation
(\ref{eq:motiondisk}),
$$
m_k  \ddot{w}_k = - \frac{2 m_k \bar{w}_k \dot{w}_k^2}{ R^2- |w_k|^2  } + \frac{( R^2- |w_k|^2  )^2}{4 R^4} \frac{\partial U_R}{\partial \bar{w}_k}, \ k=1,\dots,n,
$$
we obtain in terms of $z$ that
$$
 m_k (\ddot{z}_k  +2i \dot{z}_k -z_k ) e^{it}
$$
$$ 
=- \frac{2 m_k e^{-it} \bar{z}_k  (i z_k + \dot{z}_k )^2 e^{2it} }{ R^2- |z_k|^2  } +
 \frac{( R^2- |z_k|^2  )^2}{4 R^4} \, \frac{\partial U_R}{\partial \bar{z}_k} \,\frac{d \bar{z}_k}{d \bar{w}_k}
$$
$$
= - \frac{2 m_k \bar{z}_k   (i z_k + \dot{z}_k )^2 e^{it} }{ R^2-
|z_k|^2  } + \frac{( R^2- |z_k|^2  )^2}{4 R^4} \,
\frac{\partial U_R}{\partial \bar{z}_k} \, e^{it}.
$$
Since $z$ is a solution of (\ref{eq:motiondisk}), $m_k \neq 0$, and $e^{it} \neq 0$, the last relationship becomes
\begin{equation} \label{eq:condinvardisk1}
2i \dot{z}_k -z_k    = - \frac{2 \bar{z}_k (2i z_k  \dot{z}_k -z_k^2)
}{ R^2- |z_k|^2  } ,
\end{equation}
which is equivalent to the equation
\begin{equation} \label{eq:condinvardisk2}
2i \left[1+ \frac{2 |z_k|^2}{ R^2- |z_k|^2 }\right]\dot{z}_k =
\left[1+ \frac{2 |z_k|^2}{ R^2- |z_k|^2 }\right] z_k.
\end{equation}
Equation (\ref{eq:condinvardisk2}) holds if and only if
$$
1+\frac{2 |z_k|^2}{ R^2- |z_k|^2 }=0 \ \ {\rm or}\ \ 2i \dot{z}_k = z_k, \ \ k=1,\dots,n.
$$
The first set of conditions is equivalent to  $|z_k(t)|=0=R$, which never holds. 
The second set of conditions, which provides information about how the velocities must behave for this kind of relative equilibria, holds for $|z_k(t)|=r_k$, where  $r_k\ge 0, k=1,\dots, n$. For this reason, the particles form a relative equilibrium associated to the Killing vector field (\ref{eq:killingvector}) if they are moving along Euclidean circles centered at the origin of the coordinate system in $\mathbb{D}^2_R$. In terms of $\mathbb L_R^2$, the particles move on circles obtained by slicing the hyperbolic sphere with planes orthogonal to the vertical axis, $\mathfrak z$, of $\mathbb R_1^3$.

We can now prove the following result.

\begin{Theorem}\label{thm:existence} Consider $n$ point particles
with masses $m_1,\dots, m_n>0$, $n\ge 2$, moving in $\mathbb{D}^2_{R}$. A necessary and sufficient condition for  the function $z=(z_1,\dots, z_n)$ to be a solution of system \eqref{eq:motiondisk} that is a relative equilibrium associated to the Killing vector field $L_Z$ defined by equation \eqref{eq:killingvector} is that for all $k=1,\dots, n,$ the following equations are satisfied
\begin{equation} \label{eq:rationalsystemdisk}
\frac{R^3(R^2+r_k^2) z_k }{ 4(R^2- r_k^2)^4 }
= - \sum_{\substack{j=1\\ j\ne k }}^n \frac{m_j (r_j^2-R^2)^2(z_j-z_k)(R^2-z_k\bar{z}_j)}{(\tilde{\Theta}_{2,(k,j)}(z, \bar{z}))^{3/2}},
\end{equation}
where $r_k=|z_k|$ and 
$$
\tilde{\Theta}_{2,(k,j)}(z, \bar{z})
$$
$$
=[2(z_k \bar{z}_j+z_j \bar{z}_k)R^2-
(r_k^2+R^2)(r_j^2+R^2)]^2-(R^2-r_k^2)^2(R^2-r_j^2)^2,
$$
 $k,j\in\{1,\dots, n\}, \ k\ne j$.
\end{Theorem}
\begin{proof}
From the equations $2i \dot{z}_k = z_k, \ k=1,\dots,n$, we concluded that a necessary condition for the existence of a relative equilibrium of the aforementioned type is that the particles move along ordinary circles centered at the origin of the coordinate system in $\mathbb{D}^2_{R}$. Differentiating these conditions and using them again we obtain that
\begin{equation} \label{harmonic}
-4 \ddot{z}_k = z_k, \ k=1,\dots,n.
\end{equation}
Comparing these equalities with the equations of motion (\ref{eq:motiondisk}), we conclude that the coordinates of a relative equilibrium must satisfy the $n$ algebraic equations
\begin{equation} \label{eq:condrationalsystemdisk}
m_k z_k = -\frac{2 m_k |z_k|^2 z_k}{ R^2- |z_k|^2} - \frac{2( R^2-
|z_k|^2  )^2}{R^4} \, \frac{\partial U_R}{\partial \bar{z}_k},\ k=1,\dots,n.
\end{equation}

Substituting $r_k=|z_k|, \ k=1,\dots, n$, into the above equations, we obtain the system of $n$ equations (\ref{eq:rationalsystemdisk}), which characterize the relative equilibria given by the group $G_2$. This remark completes the proof.
\end{proof}

\begin{Definition} \label{def:eliprelequil} We call {\it elliptic relative equilibria} the solutions of system \eqref{eq:motiondisk} that satisfy the conditions \eqref{eq:rationalsystemdisk}.
\end{Definition}

\subsection{The case $n = 2$}\label{sec:case_2}

We will next prove the existence of elliptic relative equilibria for 2 particles in $\mathbb{D}^2_{R}$, both in the case when they move on the same suitable circle and in the case when they move on different suitable circles. To achieve this goal, notice first that for $n=2$ and $m_1,m_2>0$, equations \eqref{eq:rationalsystemdisk}  take the form
\begin{equation*}
\frac{R^3(R^2+r_1^2) z_1 }{4(R^2-
r_1^2)^4  } = - \frac{m_2 (r_2^2-R^2)^2(z_2-z_1)(R^2-z_1\bar{z}_2)}{[\tilde{Q}_{2,(1,2)} (z, \bar{z})]^{3/2}}
\end{equation*}
\begin{equation*}
\frac{R^3(R^2+r_2^2) z_2 }{4(R^2- r_2^2)^4  } = - \frac{m_1
(r_1^2-R^2)^2(z_1-z_2)(R^2-z_2\bar{z}_1)}{[\tilde{Q}_{2,(2,1)} (z, \bar{z})]^{3/2}}, 
\end{equation*}
where
$$
\tilde{Q}_{2,(k,j)}(z, \bar{z})= [2(z_1 \bar{z}_2+z_2 \bar{z}_1)R^2- (r_1^2+R^2)(r_2^2+R^2)]^2-(r_1^2-R^2)^2(r_2^2-R^2)^2.
$$
Some algebraic manipulations lead us to the equation
$$
\frac{(R^2+r_1^2) (R^2-r_2^2)^2  m_1
}{(R^2+r_2^2)(R^2- r_1^2)^2  m_2}= 
\frac{z_2 (z_2-z_1)(R^2-z_1\bar{z}_2)}{z_1 (z_1-z_2)(R^2-z_2\bar{z}_1)} ,
$$
and if we simplify the right hand side we obtain
\begin{equation}
\label{eq:twobodyequilibriadisk}
\frac{(R^2+r_1^2) (R^2-r_2^2)^2  m_1
}{(R^2+r_2^2)(R^2- r_1^2)^2  m_2} = - \frac{R^2 z_2 -z_1 r_2^2}{R^2 z_1 -z_2 r_1^2}.  
\end{equation}
This equation shows that there are no elliptic relative equilibria for the 2-body problem in $\mathbb{D}^2_{R}$ when one particle is fixed at the origin of the coordinate system.
Indeed, if $z_1=0$, then $r_1=0$, and the denominator of the right hand side vanishes.
If $z_2=0$, then $r_2=0$, and the above equation becomes
$$
\frac{(R^2+r_1^2)R^2m_1}{(R^2-r_1)^2m_2}=0,
$$
which has, obviously, no solutions.

Equation \eqref{eq:twobodyequilibriadisk} holds for all time, in particular 
when the particle $m_1$ reaches the real line. At that time instant, we have $z_1=\alpha\in\mathbb R$, and let us denote $z:=z_2$ and $r:=|z_2|$. Then, if we solve equation
(\ref{eq:twobodyequilibriadisk}) for $z$, we have
\begin{equation*}
z=\frac{m_1 (R^2+ \alpha^2) (R^2-r^2)^2 R^2- m_2 r^2 (R^2+r^2)(R^2- \alpha^2)^2 }{m_1 (R^2+ \alpha^2)
(R^2-r^2)^2 \alpha^2- m_2 R^2 (R^2+r^2)(R^2- \alpha^2)^2 } \, \alpha,
\end{equation*}
therefore $z$ is also a real number, so either $z=r$ or $z=-r$. 
In other words, when $m_1$ reaches the real line, $m_2$ reaches it too.
Consequently the above equation becomes
 \begin{equation}
 \frac{m_1 (R^2+ \alpha^2) (R^2-r^2)^2}{m_2 (R^2+ r^2) (R^2-\alpha^2)^2}=
 -\frac{(\pm r)[R^2-\alpha (\pm r)]}{\alpha[R^2-\alpha (\pm r)]}.
 \end{equation}
Since $\alpha, r < R$, it follows that the left hand side is positive, so the right hand side
must be positive too. Consequently $z=-r$, and thus the above equation takes the form
\begin{equation}
 \frac{m_1 (R^2+ \alpha^2) (R^2-r^2)^2 }{m_2 (R^2+ r^2) (R^2-\alpha^2)^2 }=
 \frac{r}{\alpha},
 \label{r/alfa}
 \end{equation}
which we can use to estimate the value of $r$ that gives the position of $m_2$ as a function of $m_1, m_2, R,$ and $\alpha$. For this purpose, we consider the real function $f\colon(-R,R)\to\mathbb R$,
\begin{equation}\label{eq:real-function-eliptic}
f(x)=m_1 \alpha (R^2+ \alpha^2) (R^2-x^2)^2 - m_2 x (R^2+ x^2) (R^2-\alpha^2)^2,
\end{equation}
whose zeroes give us the desired values of $r$ and, therefore, the elliptic relative
equilibria for the 2-body problem in $\mathbb{D}_R^2$. For this purpose, let us
first prove the following result.

\begin{Lemma} \label{lema:not-double-zeros-elliptic-two-body}
The function $f$ defined in \eqref{eq:real-function-eliptic} has no double roots.
\end{Lemma}
\begin{proof}
The first derivative of the function $f$ is
\begin{equation}\label{eq:derivative-elliptic-rel- equi-two-body}
f'(x)=-4m_1 \alpha x(R^2+ \alpha^2) (R^2-x^2) - m_2  (R^2+ 3x^2) (R^2-\alpha^2)^2.
\end{equation}
Since $x \neq 0$ and $-R<x<R$, the double zeroes of $f$ must satisfy the equations
$$
\begin{cases}
    -4x \, f(x) = 0\cr
    (R^2-x^2) \,f'(x) = 0, \cr
\end{cases}
$$
which are equivalent to the system
\begin{equation*}\label{eq:system-double-zeros-elliptic}
\begin{cases}
4 m_1 \alpha x (R^2+ \alpha^2) (R^2-x^2)^2 - 4 m_2 (R^2 x^2+ x^4) (R^2-\alpha^2)^2   =0\cr
-4m_1 \alpha x(R^2+ \alpha^2) (R^2-x^2)^2 - m_2  (R^2+ 3x^2) (R^2-\alpha^2)^2 (R^2-x^2)=0.\cr
\end{cases}
\end{equation*}
If we add the equations of the above system, we obtain the quartic equation
$$
x^4+6 R^2 x^2 +R^4=0,
$$
which has only non-real roots, $\pm \sqrt{-3\pm 2\sqrt{2}} \, R$, a fact that completes the proof. 
\end{proof}

We can now state and prove the following result, which characterizes the elliptic relative equilibria of the curved 2-body problem, expressed in terms of the Poincar\'e disk model, $\mathbb D_R^2$, of the hyperbolic plane.

\begin{Theorem} \label{the:elliptic-relat-equi-two-body}
Consider 2 point particles of masses $m_1, m_2>0$ moving in the Poincar\'e disk $\mathbb{D}_R^2$, whose center is the origin, $\bf 0$, of the coordinate system. Then a function $z=(z_1,z_2)$ is an elliptic relative equilibrium of system \eqref{eq:motiondisk} with $n=2$, if and only if for every circle centered at $\bf 0$ of radius $\alpha$, with $0<\alpha <R$, along which $m_1$ moves, there is a unique circle centered at $\bf 0$ of radius $r$, which satisfies $0< r <R$ and \eqref{r/alfa}, along which $m_2$ moves, such that, at every time instant, $m_1$ and $m_2$ are on some diameter of $\mathbb{D}_R^2$, with $\bf 0$ between them. Moreover, 
\begin{enumerate}
\item if  $m_2>m_1>0$ and $\alpha$ are given, then $r< \alpha$;
\item if $m_1 =m_2>0$ and $\alpha$ are given, then $r = \alpha$;
\item if $m_1>m_2>0$ and $\alpha$ are given, then $r>\alpha$.
\end{enumerate}
\end{Theorem}
\begin{proof}
If the solution $z$ of system (\ref{eq:motiondisk}) with $n=2$ is an elliptic relative equilibrium, then equations \eqref{eq:rationalsystemdisk} are satisfied, and they lead
to equation \eqref{r/alfa}, from which we can compute $r$ for given $m_1, m_2>0$ and
$\alpha$, with $0<\alpha<R$. Then, by Lemma \ref{lema:not-double-zeros-elliptic-two-body}, there is a unique $r$ as desired, so the particles move as described, and the implication follows.

To prove the converse, assume that for given $m_1, m_2>0$ and $\alpha$, with
$0<\alpha<R$, there is a unique $r$, which satisfies $0<r<R$ and \eqref{r/alfa}, such that the bodies move as described. Then the motion must be given by the function
$z=(z_1,z_2)$ with 
$$
z_1(t)=\alpha(\cos t+i\sin t),\ \ z_2(t)=-r(\cos t +i\sin t),
$$
with the relationship between $m_1, m_2, \alpha$, and $r$ given by \eqref{r/alfa}. A straightforward computation shows that this function is a solution of system (\ref{eq:motiondisk}) with $n=2$ and satisfies equations \eqref{eq:rationalsystemdisk}.

To see how the relative values of $m_1$ and $m_2$ determine the relationship between $r$ and $\alpha$, we evaluate the function $f$, defined in (\ref{eq:real-function-eliptic}), at $x=0, \alpha, R$, and obtain
$$
f(0) = \alpha m_1(R^2+ \alpha^2) R^4> 0,
$$
$$
f(\alpha) = \alpha(m_1- m_2)(R^2+ \alpha^2) (R^2-\alpha^2)^2, 
$$
$$
f(R) = -2 R^3m_2(R^2-\alpha^2)^2 < 0. 
$$
Since, by \eqref{eq:derivative-elliptic-rel- equi-two-body}, the derivative $f'$ is negative in
the entire interval  $(0,R)$, the function $f$ is strictly decreasing. The conclusion follows then from the above relationships. This remark completes the proof.
\end{proof}

\subsection{The case $n = 3$}\label{sec:case_3}

In the case of 3 particles in the Poincar\'e disk $\mathbb D_R^2$, equations  (\ref{eq:rationalsystemdisk}) become
\begin{equation} \label{eq:rationalsystemdisk-3-particles1}
\frac{R^3(R^2+r_1^2) z_1 }{ 4(R^2- r_1^2)^4 }
\end{equation}
$$
= -  \frac{m_2 (R^2-r_2^2)^2(z_2-z_1)(R^2-z_1\bar{z}_2)}{\{[2(z_1 \bar{z}_2+z_2 \bar{z}_1)R^2-(r_1^2+R^2)(r_2^2+R^2)]^2-(R^2-r_1^2)^2(R^2-r_2^2)^2\}^{3/2}}
$$
$$
 -  \frac{m_3 (R^2-r_3^2)^2(z_3-z_1)(R^2-z_1\bar{z}_3)}{\{[2(z_1 \bar{z}_3+z_3 \bar{z}_1)R^2-(r_1^2+R^2)(r_3^2+R^2)]^2-(R^2-r_1^2)^2(R^2-r_3^2)^2\}^{3/2}},
$$
\begin{equation}\label{eq:rationalsystemdisk-3-particles2}
\frac{R^3(R^2+r_2^2) z_2 }{ 4(R^2- r_2^2)^4 }
\end{equation}
$$
= -  \frac{m_1 (R^2-r_1^2)^2(z_1-z_2)(R^2-z_2\bar{z}_1)}{\{[2(z_2 \bar{z}_1+z_1 \bar{z}_2)R^2
-(r_2^2+R^2)(r_1^2+R^2)]^2-(R^2-r_2^2)^2(R^2-r_1^2)^2\}^{3/2}}
$$
$$
-\frac{m_3 (R^2-r_3^2)^2(z_3-z_2)(R^2-z_2\bar{z}_3)}{\{[2(z_2 \bar{z}_3+z_3 \bar{z}_2)R^2-(r_2^2+R^2)(r_3^2+R^2)]^2-(R^2-r_2^2)^2(R^2-r_3^2)^2\}^{3/2}},
$$ 
\begin{equation}
\frac{R^3(R^2+r_3^2) z_3 }{ 4(R^2- r_3^2)^4 }\label{eq:rationalsystemdisk-3-particles3}
\end{equation}
$$
= -  \frac{m_1 (R^2-r_1^2)^2(z_1-z_3)(R^2-z_3\bar{z}_1)}{\{[2(z_3 \bar{z}_1+z_1 \bar{z}_3)R^2-
(r_3^2+R^2)(r_1^2+R^2)]^2-(R^2-r_3^2)^2(R^2-r_1^2)^2\}^{3/2}} 
$$
 $$
 -\frac{m_2 (R^2-r_2^2)^2(z_2-z_3)(R^2-z_3\bar{z}_2)}{\{[2(z_3 \bar{z}_2+z_2 \bar{z}_3)R^2-
(r_3^2+R^2)(r_2^2+R^2)]^2-(R^2-r_3^2)^2(R^2-r_2^2)^2\}^{3/2}}.
$$

\bigskip

\subsubsection{Eulerian Solutions}

We will start the study of the case $n=3$ with the Eulerian elliptic relative equilibria in $\mathbb{D}_R^2$ for which the bodies lie on a rotating geodesic. Of course,
we can assume that this geodesic rotates around the origin of the coordinate system, so then it must be a rotating diameter of $\mathbb D_R^2$.

Assume that the particle $m_1$ is located at the center of the disk, i.e.\  $z_1=0$, and that $m_2$ reaches at some time instant the positive axis of the real line, i.e.\ $z_2=r_2=:\alpha>0$. For $m_3$, we denote $z:=z_3$ and take $|z_3|=r_3$.
Then equations \eqref{eq:rationalsystemdisk-3-particles1}, \eqref{eq:rationalsystemdisk-3-particles2}, and \eqref{eq:rationalsystemdisk-3-particles3}
become, respectively,
\begin{equation} \label{eq:rationalsystemdisk-3-particles-euler1}
0 =-  \frac{m_2(R^2-\alpha^2)^2}{\alpha^{2}} -  \frac{m_3 z (R^2-r_3^2)^2}{r_3^{3}},
\end{equation}

\begin{equation}\label{eq:rationalsystemdisk-3-particles-euler2}
\frac{R^3(R^2+\alpha^2) \alpha }{ 4(R^2- \alpha^2)^4 }
= \frac{m_1 \alpha}{[(\alpha^2+R^2)^2-(R^2-\alpha^2)^2]^{3/2}} 
\end{equation}
$$
-\frac{m_3 (R^2-r_3^2)^2(z-\alpha)(R^2-\alpha\bar{z})}{\{[2 \alpha (\bar{z}+z)R^2-
(\alpha^2+R^2)(r_3^2+R^2)]^2-(R^2-\alpha^2)^2(R^2-r_3^2)^2\}^{3/2}}, 
$$

\begin{equation}\label{eq:rationalsystemdisk-3-particles-euler3}
\frac{R^3(R^2+r_3^2) z }{ 4(R^2- r_3^2)^4 }
=\frac{m_1 z}{[
(r_3^2+R^2)^2-(R^2-r_3^2)^2]^{3/2}} 
\end{equation}
$$
- \frac{m_2 (R^2-\alpha^2)^2(\alpha-z)(R^2-\alpha z)}{\{[2 \alpha (z+\bar{z})R^2-
(r_3^2+R^2)(\alpha^2+R^2)]^2-(R^2-r_3^2)^2(R^2-\alpha^2)^2\}^{3/2}}.
$$

\medskip

The following result will show that, when one particle is fixed at the origin of the coordinate system at the center of $\mathbb D_R^2$, there is just one class of Eulerian elliptic relative equilibria, namely orbits for which the distance from the fixed body to the 2 rotating bodies is the same, and consequently those masses must be equal. In terms of the hyperbolic sphere $\mathbb L_R^2$, the configuration is an isosceles triangle that rotates around its vertical height.

\begin{Theorem}\label{theo:elliptic-relative-equilibria-3-bodies} 
Consider 3 point particles of masses $m_1, m_2, m_3>0$ moving in the Poincar\'e disk $\mathbb{D}_R^2$, whose center is the origin, $\bf 0$, of the coordinate system. Take a function $z=(z_1,z_2, z_3)$ that describes the positions of the particles, with $z_1(t)=0$ for all $t$. Then $z$ is an Eulerian elliptic relative equilibrium of system \eqref{eq:motiondisk} with $n=3$ if and only if $m_2$ and $m_3$ are at the opposite sides of the same uniformly rotating diameter of a circle of radius $\alpha$ in $\mathbb{D}_R^2$, centered at $\bf 0$, with $0<\alpha<R$, and $m_1=m_2$.
 \end{Theorem}
\begin{proof}
From equation (\ref{eq:rationalsystemdisk-3-particles-euler1}), we can conclude that 
$z$ must be a negative real number, so $z=-r_3=:-r$. This implies that \eqref{eq:rationalsystemdisk-3-particles-euler1} becomes
\begin{equation}
\frac{m_2(R^2-\alpha^2)^2}{\alpha^{2}} = \frac{m_3  (R^2-r^2)^2}{r^{2}}.
\end{equation}
Let us further consider the previous equation together with \eqref {eq:rationalsystemdisk-3-particles-euler2} and \eqref {eq:rationalsystemdisk-3-particles-euler3} in which we substitute the previous values of $z_1, r_1, z_2, r_2, z_3, r_3$. Then we obtain the new system
\begin{equation} \label{eq:rationalsystemdisk-3-particles-euler-reduced}
0 = -  \frac{m_2 (R^2-\alpha^2)^2}{\alpha^{2}} +  \frac{m_3  (R^2-r^2)^2}{r^{2}},
\end{equation}

\begin{equation}\label{eq:rationalsystemdisk-3-particles-euler-reduced2}
\frac{R^3(R^2+\alpha^2) \alpha }{ 4(R^2- \alpha^2)^4 }
=\frac{m_1}{ 2^{3/2}8 R^3 \alpha^2}
\end{equation}
$$
+\frac{m_3 (R^2-r^2)^2(r+\alpha)(R^2+\alpha r)}{\{[-4 R^2 r \alpha -
 (\alpha^2+R^2)(r^2+R^2)]^2-(R^2-\alpha^2)^2(R^2-r^2)^2\}^{3/2}}, 
 $$
 
 \begin{equation}\label{eq:rationalsystemdisk-3-particles-euler-reduced3}
 \frac{R^3(R^2+r^2) r  }{ 4(R^2- r^2)^4 }
=   \frac{m_1}{ 2^{3/2}8 R^3 r^2}
\end{equation}
$$
+\frac{m_2 (R^2-\alpha^2)^2(\alpha+r)(R^2+\alpha r)}{\{[-4 R^2 r \alpha-
(r^2+R^2)(\alpha^2+R^2)]^2-(R^2-r^2)^2(R^2-\alpha^2)^2\}^{3/2}}. 
$$

If  we multiply both sides of (\ref{eq:rationalsystemdisk-3-particles-euler-reduced2})
by $m_2(R^2-\alpha^2)^2$, and both sides of \eqref{eq:rationalsystemdisk-3-particles-euler-reduced3} by $- m_3(R^2-r^2)^2$, when we add the resulting equations
we get
\begin{equation}\label{eq:eulerian-three-2}
\frac{m_2(R^2+\alpha^2) \alpha}{(R^2-\alpha^2)^2}=\frac{m_3 (R^2+ r^2) r}{(R^2- r^2)^2}.
\end{equation}
From equations \eqref{eq:rationalsystemdisk-3-particles-euler-reduced} and (\ref{eq:eulerian-three-2}) we obtain the linear system having the masses $m_2$ and $m_3$ as unknowns,
\begin{equation}\label{eq:eulerian-three-3}
\begin{cases}
m_2 r^{2} (R^2-\alpha^2)^2 -m_3 \alpha^{2} (R^2-r^2)^2 = 0 \cr
m_2 (R^2+\alpha^2) (R^2- r^2)^2 \alpha -m_3 (R^2+ r^2) (R^2-\alpha^2)^2 r = 0, \cr
\end{cases}
\end{equation}
which has nontrivial solutions if and only if the principal determinant vanishes.
A straightforward computations shows that this condition is equivalent to
\[(R^2+ r^2) (R^2-\alpha^2)^4 r^3 -(R^2+\alpha^2) (R^2- r^2)^4 \alpha^3  =0,\]
so the principal determinant vanishes when
\begin{equation}\label{eq:eulerian-three-solution}
\frac{(R^2+ r^2) r^3}{(R^2- r^2)^4}= \frac{(R^2+\alpha^2)  \alpha^3 }{(R^2-\alpha^2)^4}.
\end{equation}
To find the values of $r$ that solve equation (\ref{eq:eulerian-three-solution}), we consider the function
\begin{equation}\label{eq:function-eulerian-three-solution}
g\colon[0,R)\to\mathbb R, \ \ \ g(x)=\frac{(R^2+ x^2) x^3}{(R^2- x^2)^4},
\end{equation}
which is strictly increasing in its domain. Therefore equation (\ref{eq:eulerian-three-solution}) holds only for $x=r=\alpha$. If we substitute these values in the first equation of system (\ref{eq:eulerian-three-3}), it follows that $m_2=m_3$. 
\end{proof}

In Theorem \ref{theo:elliptic-relative-equilibria-3-bodies}, we required that all masses are positive. Let us now consider the case when $m_1=0$. We then obtain the following result.

\begin{Proposition}\label{prop:restricted-elliptic-relative-equilibria-3-bodies} 
Consider 3 point particles of masses $m_1=0$ and $m_2= m_3>0$ moving in the Poincar\'e disk $\mathbb{D}_R^2$, whose center is the origin, $\bf 0$, of the coordinate system. Take a function $z=(z_1,z_2, z_3)$ that describes the positions of the particles, and assume that $z_1(0)=0$ and the real parts of $z_2(0)$ and $z_3(0)$ are $0$, i.e.\ $m_1$ is initially at the center and $m_2, m_3$ are initially on the horizontal diameter of $\mathbb{D}_R^2$. Then a necessary condition for the particles to form an elliptic relative equilibrium, is that $m_2$ and $m_3$ rotate on the same suitable circle, being at every time instant at the opposite sides of some diameter of that circle, and $m_1$ lies at the the center of the disk for all time.
\end{Proposition}
\begin{proof}
Without loss of generality, we can take 
$$
z_1=r_1=\bar{z}_1=:c,\ \ z_2=r_2=\bar{z}_2=:\alpha,\ \ {\rm and}\ \ z_3=\bar{z}_3=-\beta.
$$
Then some
straightforward computations show that equations (\ref{eq:rationalsystemdisk-3-particles-euler1}), (\ref{eq:rationalsystemdisk-3-particles-euler2}), and (\ref{eq:rationalsystemdisk-3-particles-euler3}) become, respectively, 
\begin{equation} \label{eq:rational-restricted-systemdisk-3-particles-euler1}
\frac{R^6(R^2+c^2)c}{(R^2-c^2)^4} = -  \frac{m_2(R^2-\alpha^2)^2}{2(\alpha-c)^2(R^2-\alpha c)^2}+
\frac{m_3  (R^2-\beta^2)^2}{2(\beta+c)^2(R^2+ \beta c)^2}, 
\end{equation}
\begin{equation}\label{eq:rational-restricted-systemdisk-3-particles-euler2}
\frac{R^6(R^2+\alpha^2)\alpha}{(R^2-\alpha^2)^4} = -\frac{m_3  (R^2-\beta^2)^2}{2(\beta+\alpha)^2(R^2+\beta \alpha)^2}, 
\end{equation}
\begin{equation}\label{eq:rational-restricted-systemdisk-3-particles-euler3}
\frac{R^6(R^2+\beta^2)\beta}{(R^2-\beta^2)^4} = -\frac{m_2 (R^2-\alpha^2)^2}{2(\beta+\alpha)^2(R^2+\beta \alpha)^2}.
\end{equation}
Since $m_2=m_3$, equations (\ref{eq:rational-restricted-systemdisk-3-particles-euler2}) and (\ref{eq:rational-restricted-systemdisk-3-particles-euler3}) lead us to the relationship
\begin{equation}\label{eq:relation-restricted-eulerian-elliptic}
\frac{(R^2+\alpha^2) \alpha}{(R^2 -\alpha^2)^2} = \frac{(R^2+\beta^2) \beta}{(R^2 -\beta^2)^2}.
\end{equation}
Consider the smooth real function
$$
h\colon[0,R) \to \mathbb{R},\ \ 
h(x)= \frac{(R^2+x^2) \, x}{(R^2 -x^2)^2}.
$$
It is easy to see that $h$ is strictly increasing, which implies that relation (\ref{eq:relation-restricted-eulerian-elliptic}) holds if and only if $\alpha=\beta$.
Therefore $m_2$ and $m_3$ must be at opposite sides of the same circle.

If we now use the fact that $\alpha=\beta$ in equation (\ref{eq:rational-restricted-systemdisk-3-particles-euler1}), its right hand side vanishes, so
 $c=0$, a remark that completes the proof. 
 \end{proof}

\subsubsection{Lagrangian Solutions}

We will next study the Lagrangian elliptic relative equilibria in $\mathbb{D}_R^2$. As
we proved in \cite{Diac}, such orbits, formed by rotating equilateral triangles, i.e.\ $r:=|z_1|=|z_2|=|z_3|$, exist only when the 3 positive masses are equal, $m:=m_1=m_2=m_3>0$. The converse is also true, and Florin Diacu gave a proof of this result for polygons with $n\ge 3$ sides in \cite{Diacu1} in the more general case of homographic orbits, which allow not only rotation but also expansion and contraction of the configuration, from which the result for relative equilibria follows as a particular case. Moreover, Pieter Tibboel recently proved that irregular polygons cannot form homographic orbits in the 2-dimensional hyperbolic space, as well as on the 2-dimensional sphere as long as the motion is not along a great circle of the sphere, \cite{Tibboel}, when the orbit is necessarily a relative equilibrium, a case in which non-equilateral triangles exist for suitable non-equal masses, as proved in \cite{Diacu1} for the curved 3-body problem. The proof we give here for the case $n=3$ follows the idea used for Theorem 5.3 in \cite{Perez} for positive curvature.

\begin{Theorem}\label{theo:lagrangian-relative-equilibria}
Assume that 3 point particles of equal masses move along a circle of radius $r$ centered at the center of the the Poincar\'e disk $\mathbb{D}_R^2$. Then a necessary and sufficient condition for the existence of an elliptic relative equilibrium is that the particles form an equilateral triangle.
\end{Theorem}
\begin{proof}
With the values $z_1=:r$, $z_2=:re^{i \theta_2}$ and $z_3=z_2=:re^{i \theta_3}$, some straightforward computations bring equations \eqref{eq:rationalsystemdisk-3-particles1}, \eqref{eq:rationalsystemdisk-3-particles2}, and
\eqref{eq:rationalsystemdisk-3-particles3}, respectively, to 
\begin{equation} \label{eq:lagrangian-3-particles1}
\frac{R^3(R^2+r^2)}{ 4 m (R^2- r^2)^6 }
\end{equation}
$$
= -  \frac{(e^{i \theta_2}-1)(R^2-r^2e^{-i \theta_2})}{\{[4 R^2 r^2 \cos \theta_2 -
(r^2+R^2)^2]^2-(R^2-r^2)^4\}^{3/2}}
$$
$$
-\frac{(e^{i \theta_3}-1)(R^2-r^2e^{-i \theta_3})}{\{[4  R^2 r^2 \cos \theta_3 -
(r^2+R^2)^2]^2-(R^2-r^2)^4\}^{3/2}},
$$
\begin{equation} \label{eq:lagrangian-3-particles2}
\frac{R^3(R^2+r^2)}{ 4 m (R^2- r^2)^6 }
\end{equation}
$$
 =-  \frac{e^{-i \theta_2 }(1-e^{i \theta_2})(R^2-r^2 e^{i \theta_2})}{\{[4  R^2  r^2 \cos \theta_2 -(r^2+R^2)^2)]^2-(R^2-r^2)^4\}^{3/2}} 
$$
$$
- \frac{e^{-i \theta_2 }(e^{i \theta_3}-e^{i \theta_2})(R^2-r^2 e^{i (\theta_2- \theta_3)})}{\{[4 R^2  r^2 \cos (\theta_2 -\theta_3) - (r^2+R^2)^2)]^2-(R^2-r^2)^4\}^{3/2}}, 
$$
\begin{equation}\label{eq:lagrangian-3-particles3}
\frac{R^3(R^2+r^2)}{ 4 m (R^2- r^2)^6 }
\end{equation}
$$
= -  \frac{e^{-i \theta_3} (1-e^{i \theta_3})(R^2-r^2 e^{i \theta_3})}{\{[4  R^2 r^2 \cos \theta_3 -(r^2+R^2)^2]^2-(R^2-r^2)^4\}^{3/2}} 
$$
$$
- \frac{e^{-i \theta_3}(e^{i \theta_2}-e^{i \theta_3})(R^2-r^2 e^{i (\theta_3-\theta_2)})}{[[4  R^2 r^2 \cos (\theta_3 -\theta_2) -
(r^2+R^2)^2)]^2-(R^2-r^2)^4]^{3/2}}. 
$$
Adding equations \eqref{eq:lagrangian-3-particles1} and \eqref{eq:lagrangian-3-particles2}, subtracting equation \eqref{eq:lagrangian-3-particles3} from the
sum, and separating the real and imaginary parts of the resulting equation, we obtain
\begin{equation}\label{eq:trascen-algebra2}
\frac{1-\cos \, \theta_2}{D_{12} }+ \frac{1-\cos (\theta_3 -\theta_2)}{D_{23}} = \frac{2(1-\cos \,\theta_3)}{D_{13}} 
\end{equation}
\begin{equation}\label{eq:trascen-algebra21}
\frac{\sin \, \theta_2}{D_{12} }+ \frac{\sin (\theta_3 -\theta_2)}{D_{23}} = 0, 
\end{equation}
where
$$
D_{12} = 8^{3/2} R^3 r^3 [1- \cos \, \theta_2]^{3/2}[R^4+r^4- 2 R^2 r^2 \cos \, \theta_2]^{3/2}, 
$$
$$
D_{23} = 8^{3/2} R^3 r^3 [1- \cos (\theta_3- \theta_2)]^{3/2}[R^4+r^4- 2 R^2  r^2 \cos (\theta_3- \theta_2)]^{3/2}, 
$$
$$
D_{13} = 8^{3/2} R^3 r^3 [1- \cos \, \theta_3]^{3/2}[R^4+r^4- 2  R^2  r^2 \cos \, \theta_3]^{3/2}. 
$$
It is easy to check that $\theta_2=\frac{2 \pi}{3}$ and $\theta_3=\frac{4 \pi}{3},$ satisfy equations (\ref{eq:trascen-algebra2}) and therefore the configuration corresponds to an equilateral triangle.

Using standard trigonometry, equation (\ref{eq:trascen-algebra21}) becomes
\[ \left[\frac{\sin^2 (\frac{\theta_3- \theta_2}{2})}{\sin^2 (\frac{\theta_2}{2})} \right]^2
\left[\frac{(r^2-R^2)^2+R^2r^2 \sin^2 (\frac{\theta_3-
\theta_2}{2})}{(r^2-R^2)^2+R^2r^2 \sin^2
(\frac{\theta_2}{2})}\right]^3 = \frac{1-\sin^2 (\frac{\theta_3-
\theta_2}{2})}{1-\sin^2 (\frac{\theta_2}{2})}.
\]
Renaming the variables as $u=\sin^2 (\frac{\theta_3- \theta_2}{2})$
and $v=\sin^2 (\frac{\theta_2}{2})$,  the above equation takes the form
\[u^2 (1-v)[(r^2-R^2)^2+R^2 r^2 u]^3 = v^2 (1-u)[(r^2-R^2)^2+R^2 r^2 v]^3. \]
This  real  equation holds only when $u=v$, that is, $\sin^2
(\frac{\theta_3- \theta_2}{2})=\sin^2 (\frac{\theta_2}{2})$, or
equivalently, $\theta_3= 2 \theta_2$. If we substitute these values
in equation of (\ref{eq:trascen-algebra2}), we obtain
\[\frac{1- \cos \, \theta_2}{1- \cos \, 2\theta_2} = \frac{(R^4+r^4- 2  R^2 r^2 \cos \, 2\theta_2)^3}{(R^4+r^4- 2 R^2 r^2 \cos \, \theta_2)^3}.\]
Taking $w=\cos \, \theta_2$ and $s=\cos \, 2\theta_2$,  we are
led to
\[(1-w)(R^4+r^4- 2R^2 r^2 s)^3 = (1-s)(R^4+r^4- 2R^2 r^2 w)^3,  \]
which has real solutions only for $w=s$, i.e.\ when $\cos \, \theta_2=
\cos \, 2\theta_2$, which yields
$$\theta_2= 0, \, \frac{2 \pi}{3}, \, \frac{4\pi}{3}, \, 2 \pi. $$

To avoid singular configurations, we must take $\theta_2=
\frac{2 \pi}{3}, \, \frac{4\pi}{3}$, which correspond to an
equilateral triangle positioned in the suitable circle of radius $r$.
This remark completes the proof. 
\end{proof}

\section{The curved $n$-body problem in $\mathbb{H}^2_R$}
\label{sec:Klein_model}
Our attempts to study the relative equilibria corresponding to the subgroups $G_1$ and $G_3$ in the Poincar\'e disk $\mathbb D_R^2$ led to insurmountable computations. Therefore we chose to move the problem to the Poincar\'e's upper-half-plane model, $\mathbb{H}^2_R$, to see if the those equilibria would be easier to approach. In the
remaining part of the paper, we will obtain the equations of motion in the Poincar\'e upper half plane and analyze the hyperbolic and parabolic relative equilibria for $n=2$ and $n=3$.

\subsection{Equations of motion in $\mathbb{H}^2_R$}

We will next obtain the equations of motion of the curved $n$-body problem in Poincar\'e's upper half plane model, $\mathbb{H}^2_R$, with the help of a global isometric fractional linear transformation defined by
\begin{equation} \label{eq:transf-disk-plane}
z\colon\mathbb{H}^2_R \to \mathbb{D}^2_R, \ \ z= z(w)= \frac{-R w+i R^2}{w+iR},
\end{equation}
where $w$ is the complex variable in the upper half plane
\[\mathbb{H}^2_R = \{w \in \mathbb{C} \, | \, {\rm Im} \,(w) >0 \}. \]
This transformation has the inverse
$$
w\colon{\mathbb D_R^2}\to{\mathbb H_R^2},\ \ w=w(z)=\frac{iR(R-z)}{R+z}.
$$
Since
\[ dz= \frac{-2R^2i}{(w+iR)^2} dw \quad {\rm and } \quad
d\bar{z}= \frac{2R^2i}{(\bar{w}-iR)^2}  d\bar{w},\] 
the metric (\ref{met-c}) of the disk $\mathbb{D}^2_R$ is transformed into the metric
\begin{equation}\label{met-k}
   -ds^2= \frac{4 R^2}{(w-\bar{w})^2}dw d\bar{w}.
\end{equation}
Then $\mathbb{H}^2_R$ endowed with the metric (\ref{met-k}) is called the Poincar\'e upper half plane model of hyperbolic geometry, for which the conformal factor is given by
\[ \mu (w, \bar{w}) = -\frac{4 R^2 }{(w-\bar{w})^2}. \]
In terms of the metric (\ref{met-k}), the geodesics are either half circles orthogonal to the real axis $(y=0)$ or half lines perpendicular to it.  All these curves have infinite length.

Applying transformation (\ref{eq:transf-disk-plane}) to equation (\ref{eq:potesf}), we obtain the new potential in the coordinates $(w, \bar{w})$, given by
\begin{equation}\label{eq:potesfklein}
V_R (w, \bar{w})= \frac{1}{R} \sum_{1 \leq k < j \leq n}^n
m_k m_j \frac{(\bar{w}_k+w_k)(\bar{w}_j+w_j)-2(|w_k|^2+|w_j|^2)}{[\Theta_{3,(k,j)}(w, \bar{w})]^{1/2}},
\end{equation}
where
\begin{equation}\label{eq:singularklein}
\Theta_{3,(k,j)}(w, \bar{w})
\end{equation}
$$
=[(\bar{w}_k+w_k)(\bar{w}_j+w_j)-2(|w_k|^2+|w_j|^2)]^2-(\bar{w}_k-w_k)^2(\bar{w}_j-w_j)^2.
$$
Notice that with the values $w_k= x_k+iy_k $ and $w_j= x_j+iy_j$, we obtain 
\begin{equation}\label{positive}
\Theta_{3,(k,j)}
\end{equation}
$$
= 4(x_k-x_j)^2[(x_k-x_j)^2+ 2(y_k^2+y_j^2)]+4(y_k^2-y_j^2)^2,
$$
which is always positive, except when it takes the value zero as collisions take place, i.e.\ when $x_k=x_j$ and $y_k=y_j$ for at least 2 indices $k,j\in\{1,\dots,n\}, k\ne j$.


A result analogous to Lemma \ref{lema:gral-mechanicalhyp} in the Poincar\'e disk, $\mathbb D_R^2$, is the following property, valid in the Poincar\'e upper half plane, $\mathbb H_R^2$.

\begin{Lemma}\label{lema:gral-mechanicalhypklein}  Let
\[ \displaystyle L(w,\dot{w})= \frac{1}{2}\sum_{k=1}^n m_k \, \mu(w_k, \bar{w}_k) \, |\dot{w}_k|^2 + V_R(w, \bar{w})  \]
be the Lagrangian for the given problem in $\mathbb{H}^2_R$. Then the  solution curves of the corresponding Euler-Lagrange equations satisfy the system of $n$ second order differential equations, $ k=1,\dots,n$,
\begin{equation} \label{eq:motionkleingeo}
  m_k  \ddot{w}_k - \frac{2 m_k  \dot{w}_k^2}{ w_k-\bar{w}_k} =\frac{2}{\mu(w_k, \bar{w}_k)} \,
  \frac{\partial V_R}{\partial \bar{w}_k} =-\frac{( w_k- \bar{w}_k )^2}{2 R^2} \,
  \frac{\partial V_R}{\partial \bar{w}_k},
\end{equation}
where
\begin{equation}\label{eq:gradmet2}
 \frac{\partial V_R}{\partial \bar{w}_k}= \sum_{\substack{j=1\\ j\ne k}}^n \frac{4 m_k m_j R (\bar{w}_k-w_k)(\bar{w}_j-w_j)^2(w_k-w_j)(\bar{w}_j-w_k)}{(\Theta_{3,(k,j)}(w, \bar{w}))^{3/2}}
\end{equation}
and $\Theta_{3,(k,j)}(w, \bar{w})$ is defined by \eqref{eq:singularklein}.
\end{Lemma}
\begin{proof}
From the definition (\ref{eq:transf-disk-plane}) of the linear fractional transformation
$z$, we obtain in terms of the coordinates that
$$
\dot{z}_k = \frac{-2 R^2 i}{(w_k+iR)^2} \dot{w}_k,
$$
$$
\ddot{z}_k = \frac{4 R^2 i}{(w_k+iR)^3} \dot{w}_k^2-\frac{2 R^2 i}{(w_k+iR)^2} \ddot{w}_k,
$$
$$
\frac{\partial \bar{w}_k}{\partial \bar{z}_k}= \frac{(\bar{w}_k -i R)^2}{2R^2 i}.
$$
Therefore the equations of motion (\ref{eq:motiondiskgeo}) of the curved $n$-body problem in the Poincar\'e disk get transformed into system (\ref{eq:motionkleingeo}) in the Poincar\'e upper half plane. This remark completes the proof. 
\end{proof}

\subsection{Integrals of motion in $\mathbb H_R^2$}

We will further compute the 4 integrals of system \eqref{eq:motionkleingeo} by 
transforming the 4 integrals of system \eqref{eq:motiondisk} via the transformation
\eqref{eq:transf-disk-plane}.

The energy integral is given by
\begin{equation}\label{energyH}
\frac{1}{2}\sum_{k=1}^n m_k \, \mu(w_k, \bar{w}_k) \, |\dot{w}_k|^2 - V_R(w, \bar{w})=h,
\end{equation}
where $V_R$ is given by \eqref{eq:potesfklein} and $h$ is the energy constant, the same that occurs in \eqref{Deng}, which was used to obtain \eqref{energyH}.

With the transformation \eqref{eq:transf-disk-plane}, the 3 integrals of the total angular momentum, \eqref{Dang1}, \eqref{Dang2}, and \eqref{Dang3}, become 
\begin{equation}\label{Hang1}
\sum_{k=1}^n\frac{m_kR[|w_k|^2+R^2-iR(w_k-\bar{w}_k)]^2}{2(w_k-\bar{w}_k)^2(w_k+iR)^2(\bar{w}_k-iR)^2}[\dot{\bar{w}}_kw_k^2-\dot{w}_k\bar{w}_k^2-R^2(\dot{\bar{w}}_k-\dot{w}_k)]=c_1,
\end{equation}
\begin{equation}\label{Hang2}
\sum_{k=1}^n\frac{m_kR^2[|w_k|^2+R^2-iR(w_k-\bar{w}_k)]^2}{(w_k-\bar{w}_k)^2(w_k+iR)^2(\bar{w}_k-iR)^2}(\dot{w}_k\bar{w}_k+\dot{\bar{w}}_kw_k)=c_2,
\end{equation}
\begin{equation}\label{Hang3}
\sum_{k=1}^n\frac{m_kR[|w_k|^2+R^2-iR(w_k-\bar{w}_k)]^2}{2(w_k-\bar{w}_k)^2(w_k+iR)^2(\bar{w}_k-iR)^2}[\dot{\bar{w}}_kw_k^2+\dot{w}_k\bar{w}_k^2+R^2(\dot{\bar{w}}_k+\dot{w}_k)]=c_3,
\end{equation}
which are the integrals of the total angular momentum for system \eqref{eq:motionkleingeo}. A straightforward computation confirms that the left hand sides of equations \eqref{Hang1}, \eqref{Hang2}, and \eqref{Hang3} are real functions. Moreover, the constants $c_1, c_2, c_3\in\mathbb R$ are the same that occur in equations
\eqref{Dang1}, \eqref{Dang2}, and \eqref{Dang3}.

\subsection{Relative equilibria in $\mathbb{H}^2_R$}
\label{subsec:hyperbolic-parabolic-relative eq}

In this subsection we give conditions for the existence of hyperbolic and parabolic relative equilibria in the negative curved problem by using the Poincar\'e upper half plane model, $\mathbb{H}^2_{R}$. The definitions of these concepts are given in the same geometric terms we used in Definitions \ref{def:equilibria} and \ref{def:eliprelequil}.

Let
\[ {\rm SL}(2,\mathbb{R}) = \{ A \in {\rm GL}(2,\mathbb{R}) \, | \, \det A=1 \},  \]
be the {\it special linear real 2-dimensional group}, which is a 3-dimensional simply connected, smooth real manifold. It is well known (see, e.g., \cite{Dub})  that the {\it
group of proper isometries} of $\mathbb{H}^2_R$ is the projective quotient
group $\displaystyle {\rm SL}(2,\mathbb{R})/\{\pm I \}$. Every class
\[ A = \left(\begin{array}{cc}
    a  &  b     \\
    c  &  d   \\
    \end{array}\right)  \in {\rm SL}(2,\mathbb{R})/\{\pm I \} \]
has also associated  a unique  M\"obius
transformation  $ f_A : \mathbb{H}^2_R \to \mathbb{H}^2_R$, where
\[f_A (z) = \frac{a z +b}{c z + d}, \]
for which it is easy to see that $f_{-A} (z)= f_A (z)$.
The {\it Lie algebra} of ${\rm SL}(2,\mathbb{R})$ is the 3-dimensional real
linear space
\[ {\mathfrak sl}(2,\mathbb{R}) = \{ X \in {\rm M}(2,\mathbb{R}) \, | \, \,  {\rm trace\  \! X}=0  \},  \]
spanned by the following suitable set of Killing vector fields,
\[ \left\{
X_1 = \frac{1}{2} \left( \begin{array}{ccc}
    1 & 0 \\
     0 & -1  \\
    \end{array}\right),  \quad
X_2 = \left(\begin{array}{cc}
     0 & 1 \\
    0 & 0  \\
    \end{array}\right), \quad
X_3 = \left(\begin{array}{cc}
    0 & 1 \\
    -1 & 0  \\
    \end{array}\right)
\right\} \]
As in section \ref{subsec:relativesphere}, we consider the
exponential map of matrices,
\[ \exp:  {\mathfrak sl}(2,\mathbb{R}) \to {\rm SL}(2,\mathbb{R}), \]
applied to  the one-parameter  additive subgroups (straight lines) $\{ t X_1 \}$, $ \{ t X_2 \}$,  and $\{ t X_3 \}$, to obtain the following one-parameter subgroups of the Lie group ${\rm SL}(2,\mathbb{R}) $:
\begin{enumerate}
\item  The isometric dilatation subgroup
\[ \phi_1(t)=\exp (t X_1)=
\left(\begin{array}{cc}
    e^{t/2} & 0 \\
    0  & e^{-t/2}  \\
    \end{array}\right),
\]
which defines the one-parameter family of acting M\"obius transformations
\begin{equation}\label{eq:Mobius-Klein-1}
 f_{1} (w,t) = e^{t}  w(t);
 \end{equation}

\item  The isometric shift subgroup
\[ \phi_2 (t)=\exp (t X_2)=
\left(\begin{array}{cc}
    1 & t \\
    0  & 1 \\
    \end{array}\right),
\]
which defines the one-parameter family of acting M\"obius transformations
\begin{equation}\label{eq:Mobius-Klein-2}
f_{2} (w,t) =  w(t) + t;
\end{equation}

\item  The isometric rotation subgroup
\[ \phi_3(t)=\exp (t X_3)=
\left(\begin{array}{cc}
    \cos t & \sin t \\
    -\sin t  & \cos t  \\
    \end{array}\right),
\]
which defines the one-parameter family of acting M\"obius transformations
\begin{equation}\label{eq:Mobius-Klein-3}
f_{3} (w,t) = \frac{ (\cos t) \, w(t) + \sin t }{ (-\sin t) \, w(t) + \cos t}.
\end{equation}
\end{enumerate}

We proved in subsection \ref{subsec:elliptic-relative eq} the existence of elliptic relative equilibria for the initial problem by using the proprieties of the Poincar\'e model $\mathbb{D}^2_{R}$, i.e.\ we showed that for elliptic relative equilibria, each particle moves along a circle centered at the origin of the coordinate system.

From the {\it Theorem of the invariance of the domain}, \cite{Guille}, the isometry (\ref{eq:transf-disk-plane}) carries the interior of $\mathbb{D}^2_{R}$ into the interior of $\mathbb{H}^2_{R}$. Therefore simple closed curves contained in $\mathbb{D}^2_{R}$ are taken to simple closed curves contained in $\mathbb{H}^2_{R}$ (see \cite{Guille} for more details). For the circle $z(t) = z_0 e^{it}$ in the Poincar\'e disk, $\mathbb D_R^2$, the corresponding curve,
\[ w(t)= \frac{iR(R-z_0 e^{it})}{R+ z_0 e^{it}}, \]
in the Poincar\'e upper half plane, $\mathbb H_R^2$, satisfies $w(0)=w(2 \pi)$ and  must therefore belong to the class of the M\"obius transformations (\ref{eq:Mobius-Klein-3}), which thus corresponds to elliptic relative equilibria. Since we already studied those orbits in the Poincar\'e disk, $\mathbb D_R^2$, we don't need to further analyze them here, the more so since the M\"obius transformations (\ref{eq:Mobius-Klein-3}) lead to complicated computations in $\mathbb H_R^2$. However, the M\"obius transformations (\ref{eq:Mobius-Klein-1}) and (\ref{eq:Mobius-Klein-2}), corresponding to the Killing vector fields $X_1$ and $X_2$, respectively, are simpler in $\mathbb H_R^2$ than the analogue transformations in $\mathbb D_R^2$, so they will be the object of our further analysis. The former will lead us to hyperbolic relative equilibria and the latter to parabolic relative equilibria.

 \subsection{Hyperbolic relative equilibria}
\label{subsec:hyperbolic-relative eq}

We will further study the relative equilibria associated to the subgroup
(\ref{eq:Mobius-Klein-1}), which defines the one-parameter family of acting M\"obius
transformations
\[ f_{1} (w,t) = e^{t} w (t) \]
in the Poincar\'e upper half plane, $\mathbb{H}^2_R$. Let $\xi=(\xi_1,\dots, \xi_n)$, with $\xi_k (t) = e^{t} w_k (t),\ k=1,\dots, n$, be the action orbit for a solution $w=(w_1,\dots, w_n)$ of system (\ref{eq:motionkleingeo}).
Then,
$$
\dot{\xi}_k = (w+\dot{w}) e^{t},  \ \ \
\ddot{\xi}_k = (\ddot{w}+ 2\dot{w} + w ) e^{t},  \ \ k=1,\dots,n,
$$
and therefore the curve $\xi$ is  also a solution of the equations of motion (\ref{eq:motionkleingeo}) if and only if
$$
 m_k (\ddot{w}_k  +2 \dot{w}_k +w_k ) e^{t}  = \frac{2 m_k e^{2t}  (w_k + \dot{w}_k )^2}{ e^t w_k- e^t\bar{w}_k} -
 \frac{(e^t w_k- e^t\bar{w}_k)^2}{2 R^2} \, \frac{\partial V_R}{\partial \bar{w}_k} \,\frac{d \bar{w}_k}{d \bar{\xi}_k}
$$
$$
= \frac{2 m_k (w_k + 2 w_k \dot{w}_k + \dot{w}_k^2) e^{t}}{ w_k- \bar{w}_k } - \frac{(w_k-\bar{w}_k )^2}{2 R^2} \,
\frac{\partial V_R}{\partial \bar{w}_k} \, e^{t}, \ k=1,\dots, n.  
$$
Since $w$ is a solution of (\ref{eq:motionkleingeo}), $\displaystyle \frac{d \bar{w}_k}{d \bar{\xi}_k}=e^{-t} \neq 0$, and $m_k \neq 0,\ k=1,\dots,n$, the last relationships become
\begin{equation} \label{eq:condiklein1}
2\dot{w}_k -\frac{4 w_k \dot{w}_k}{w_k- \bar{w}_k}   = - \frac{2 \bar{w}_k}{ w_k- \bar{w}_k  }-w_k, \ k=1,\dots, n.
\end{equation}
So if we fix a body $m_k$, the above condition holds if and only if
\[ 1 - \frac{2 w_k}{ w_k- \bar{w}_k }=0 \]
or
\begin{equation}\label{eq:principalcondklein1}
  2 \dot{w}_k = - w_k.
\end{equation}

The first condition in equation is equivalent to  $w_k + \bar{w}_k =0$ and corresponds to the geodesic given by the vertical half line that forms the
imaginary axis. Via the linear fractional transformation (\ref{eq:transf-disk-plane}), this geodesic corresponds  to the horizontal geodesic, ${\rm Im} (z)=0$, of the Poincar\'e disk $\mathbb{D}^2_R$.

The second condition, (\ref{eq:principalcondklein1}), holds for the body $m_k$ when $w_k$ is a function of the form $w_k(t)= w_k(0)e^{- t/2} $, where $w_k(0)$ is some initial condition with ${\rm Im} (w_k(0)) >0$. Therefore the particle $m_k$ moves along a half line through the origin of the coordinate system, which is a point at infinity that is reached when $t\to\infty$. As $t\to -\infty$, $m_k$ goes to infinity. It is instructive to notice that such integral curves of equation (\ref{eq:principalcondklein1}) correspond in the Poincar\'e disk $\mathbb{D}^2_R$ to the parametric curves
$$
z_k(t)= \frac{-R w_k(0) e^{- t/2} + iR^2}{ w_k(0) e^{- t/2} + iR}=
\frac{iR^2 e^{t/2}-R w_k(0)}{iRe^{t/2}+ w_k(0)}, 
$$ 
which are equidistant from the horizontal geodesic diameter, all of them starting at
the point $(-R,0)$, as $t \to -\infty$, and ending in the point $(R,0)$, as $t\to\infty$. 
When Re$[w_k(0)]=0$, this curve becomes the horizontal diameter, which is
a geodesic. The above non-geodesic curves together with the geodesic horizontal diameter foliate the Poincar\'e disk.

On the hyperbolic sphere, $\mathbb{L}_R^2$,  the above non-geodesic curves are given by the equations
\begin{equation}\label{eq:hyperb-in-L2}
\mathfrak x=\alpha, \quad \mathfrak y=r \sinh t, \quad \mathfrak z= r \cosh t,
\end{equation}
where $\alpha\ne 0$ is a constant and $r=\sqrt{R^2+\alpha^2}$. When $\alpha=0$,
we recover the geodesic. When these curves are taken to the Poincar\'e disk via the stereographic projection (\ref{ste-eq}), we obtain the complex curves
\begin{equation} \label{eq:project-hyperb-in-D2}
z(t)=u(t)+iv(t)=\frac{2Rr (e^{2t}-1)}{Re^{t} +2r(e^{2t}+1)} + \frac{R
\alpha e^{t} \, i}{Re^{t} +2r(e^{2t}+1)},
\end{equation}
which also start at the point $(-R,0)$, as $t\to -\infty$, and end at the point $(R,0)$, as $t \to \infty$. In other words, we have obtained two families of equidistant curves with the same initial and final directions as the geodesic diameter. 

So for the particles $m_1,\dots, m_n$ to form a relative equilibrium in $\mathbb H_R^2$ associated to the Killing vector field $X_1$, they have to move along the upper half lines converging to the origin of the coordinate system in $\mathbb{H}^2_R$ as $t \to \infty$. These half lines are not geodesics, except in the case of the vertical half line. Moreover, each non-geodesic vertical half line is equidistant from the geodesic vertical half line, the distance being larger when the angle between the non-geodesic half line and the geodesic vertical half line is larger. The sizes of these angles range between $0$ and $\pi/2$. Notice also that in the case of the geodesic vertical half line, equation
\eqref{eq:principalcondklein1} is satisfied as well, a fact that we will use in the proof of
the following result.

\begin{Theorem}\label{thm:existenceklein1} Consider $n$ point particles with masses $m_1, \dots, m_n>0$, $n\ge 2$, moving in $\mathbb{H}^2_{R}$. A necessary and sufficient  condition for the function $w=(w_1, \dots, w_n)$ to be a solution of system
\eqref{eq:motionkleingeo} and, at the same time, a relative equilibrium associated to the Killing vector field $X_1$ defined by equation \eqref{eq:Mobius-Klein-1} is that, for every $k=1,\dots,n,$ the coordinates satisfy the conditions
\begin{equation} \label{eq:condrationalsystemklein1a}
\frac{R (w_k+\bar{w}_k) \, w_k}{8(w_k- \bar{w}_k)^4} =  \sum_{\substack{j=1\\ j\ne k}}^n
\frac{m_j(w_j-\bar{w}_j)^2(w_k-w_j)(\bar{w}_j-w_k)}{[\Theta_{3,(k,j)}(w,
\bar{w})]^{3/2}},
\end{equation}
where
\begin{equation}
\Theta_{3,(k,j)}(w,
\bar{w})
\end{equation}
$$
=[(\bar{w}_k+w_k)(\bar{w}_j+w_j)-2(|w_k|^2+|w_j|^2)]^2-(\bar{w}_k-w_k)^2(\bar{w}_j-w_j)^2,
$$
$k,j\in\{1,\dots, n\}, \ k\ne j$.
\end{Theorem}
\begin{proof}
We showed previously that for relative equilibria of the aforementioned type the bodies move along straight half lines converging to the origin of the coordinate system and must therefore satisfy equation (\ref{eq:principalcondklein1}), which implies that
\begin{equation} \label{harmonicklein1}
4 \ddot{w}_k = w_k, \ k=1,\dots,n.
\end{equation}
Using this equation together with equation \eqref{eq:principalcondklein1}, we can 
conclude from the equations of motion \eqref{eq:motionkleingeo} that
\begin{equation} \label{eq:condrationalsystemklein1}
\frac{R^2m_k (w_k+\bar{w}_k) \, w_k}{2(w_k- \bar{w}_k)^3} = \frac{\partial V_R}{\partial \bar{w}_k}, \ k=1,\dots, n.
\end{equation}
Using \eqref{eq:gradmet2}, we obtain the relationships given in \eqref{eq:condrationalsystemklein1a}. This remark completes the proof.
\end{proof}

\begin{Definition} \label{def:eliprelequil} We will call {\it hyperbolic relative equilibria} the solutions of system \eqref{eq:motionkleingeo}  in $\mathbb{H}^2_{R}$ that satisfy equations \eqref{eq:condrationalsystemklein1a}.
\end{Definition}

We remark that equation (\ref{eq:principalcondklein1}) also gives the condition the velocity of the particle $m_k$ must satisfy in order to produce a hyperbolic relative equilibrium.

\subsubsection{The case $n = 2$} \label{sec:casehyp_2}

We will further provide a description of the hyperbolic relative equilibria for 2
interacting particles in $\mathbb{H}^2_{R}$. For this, we observe that
for the particles of masses $m_1$ and $m_2$, the equations 
(\ref{eq:condrationalsystemklein1a}), which characterize hyperbolic relative equilibria become
\begin{equation}
\label{eq:tworationalsystemklein}
\frac{R (w_1+\bar{w}_1) \, w_1}{8(w_1- \bar{w}_1)^4} =
\frac{m_2(\bar{w}_2-w_2)^2(w_1-w_2)(\bar{w}_2-w_1)}{[\Theta_{3,(1,2)}(w, \bar{w})]^{3/2}},
\end{equation}
\begin{equation}
\frac{R (w_2+\bar{w}_2) \, w_2}{8(w_2- \bar{w}_2)^4} =
\frac{m_1(\bar{w}_1-w_1)^2(w_2-w_1)(\bar{w}_1-w_2)}{[\Theta_{3,(2,1)}(w, \bar{w})]^{3/2}}, 
\end{equation}
where 
$$
\Theta_{3,(1,2)}(w,\bar{w})=\Theta_{3,(2,1)}(w, \bar{w})=4(w_2-w_1)(\bar{w}_1-\bar{w}_2)(\bar{w}_1-w_2)(\bar{w}_2-w_1).
$$
Straightforward computations lead to the equation
\begin{equation}
\frac{(w_1+\bar{w}_1)w_1}{(w_2+\bar{w}_2)w_2}=-\frac{m_2(\bar{w}_1-w_1)^2(\bar{w}_2-w_1)}{m_1(\bar{w}_2-w_2)^2(\bar{w}_1-w_2)},
\end{equation}
provided that $w_2+\bar{w}_2\ne 0$, which is equivalent with the equation
\begin{equation}
\label{eq:twopolinomialklein}
m_1 (w_2-\bar{w}_2)^2(\bar{w}_1-w_2)(w_1+\bar{w}_1)w_1
\end{equation}
$$
 +  m_2(w_1-\bar{w}_1)^2(\bar{w}_2-w_1)(\bar{w}_2+w_2)w_2 =0,
$$
subject to the restriction $w_2+\bar{w}_2\ne 0$.

Let us first prove a negative result about hyperbolic relative equilibria in the curved 2-body problem in $\mathbb H_R^2$. In general terms, unrelated to any model of hyperbolic geometry, this result states that, on one hand, 2 particles cannot follow each other along a geodesic and maintain a constant distance between each other and, on the other hand, one particle cannot move along the geodesic, while the other particle moves along a non-geodesic curve equidistant from that geodesic, such that the bodies maintain all the time the same distance between each other.

\begin{Proposition} \label{Prop:no-hyperbolic-equilibria-1}
Consider 2 point particles of masses $m_1, m_2>0$ moving  in $\mathbb H_R^2$. Then there are no hyperbolic relative equilibria as solutions of system \eqref{eq:motionkleingeo} with $n=2$ for which both particles move along the geodesic vertical half line on the imaginary axis or for which one particle moves along the geodesic vertical half line and the other particle moves along a non-geodesic half line converging to the origin of the coordinate system.
\end{Proposition}
\begin{proof}
As previously seen, if the particle of mass $m_1$ moves along geodesic vertical half line, then $w_1+\bar{w}_1=0$. Therefore equation (\ref{eq:twopolinomialklein}) becomes
\[
m_2(\bar{w}_1-w_1)^2(\bar{w}_2-w_1)(\bar{w}_2+w_2) \, w_2 =0,
\]
which is impossible because none of the above factors vanishes in either of the two scenarios given in the statement. This remark completes the proof. 
\end{proof}

The previous result showed that 2 bodies that form a hyperbolic relative equilibrium cannot move along the same geodesic. It is then natural to ask whether they could
move along the same non-geodesic curve that is equidistant from a given geodesic and maintain all the time the same distance from each other. As we show in the following result, expressed in terms of $\mathbb H_R^2$, the answer is also negative.

\begin{Proposition} \label{Prop:no-hyperbolic-equilibria-2}
Consider 2 point particles of masses $m_1, m_2>0$ moving in $\mathbb H_R^2$. Then there are no hyperbolic relative equilibria as solutions of system \eqref{eq:motionkleingeo} with $n=2$ for which both particles move along the same half line converging to the origin of the coordinate system.
\end{Proposition}
\begin{proof}
Recall that for hyperbolic relative equilibria the expression of the orbit of the particle $m_k$ is $w_k(t)= w_k(0)e^{- t/2}$, where ${\rm Im} (w_k(0))>0$. Then, if 2 particles belong to the same half line and form a hyperbolic relative equilibrium, their corresponding solutions have the form
\begin{equation}\label{eq:type-of sol-hyp-12}
w_1(t) = w_1(0) e^{-t/2} \,\, \,  {\rm and} \, \, \,  w_2(t)= \alpha w_1(t)= \alpha w_1(0) e^{-t/2},
\end{equation}
for some $\alpha>0$. If we substitute relations (\ref{eq:type-of sol-hyp-12}) into
equation (\ref{eq:twopolinomialklein}), a straightforward
computation gives us the condition
\[ (\alpha m_1 + m_2)w_1(0) = (m_1 + \alpha m_2)\bar{w}_1(0).
\]
If Re$(w_1(0))\ne 0$, i.e.\ when the particles are on the same non-geodesic half line, then there is no $\alpha$ that can satisfy the equation. This remark completes the proof. If Re$(w_1(0))=0$, then the bodies are on the geodesic vertical half line, and only $\alpha =-1$ satisfies the equation. But since $\alpha$ must be positive, the equation is not satisfied either, so we found another proof for the first statement of Proposition \ref{Prop:no-hyperbolic-equilibria-1}.
\end{proof}

From Propositions  \ref{Prop:no-hyperbolic-equilibria-1} and \ref{Prop:no-hyperbolic-equilibria-2} we learned that, in general terms, independent of the hyperbolic model used, 2 particles cannot form a hyperbolic relative equilibrium if they move along the same non-geodesic curve equidistantly places from a geodesic, along the same geodesic, or if one particle moves on a geodesic and the other particle on a non-geodesic curve equidistant to the geodesic. It is then natural to check the last possibility, whether there exist hyperbolic relative equilibria with one body moving along a non-geodesic curve and the other along another non-geodesic curve, both equidistant from a given geodesic. The answer is positive and the motion takes place when the distances from the geodesic to the 2 non-geodesic curves satisfy a certain relationship that depends on the values of the masses.

To write these conditions in the language of $\mathbb H_R^2$, recall that the slope $\beta$ of any straight line in the complex plane is defined by the formula
$$
i \beta =\frac{w-\bar{w}}{w+\bar{w}}.
$$
 Without loss of generality, we choose the initial conditions such that the heights satisfy $w_1(0)-\bar{w}_1(0)= 2i y_1$ and $w_2(0)-\bar{w}_2(0)= 2i y_2$. Then, in terms of the slopes $\beta_1$ and $\beta_2$ of the straight lines, equation (\ref{eq:twopolinomialklein}) becomes
\begin{equation}
\label{eq:twopolinomialkleinslopes}
 m_1  \, \beta_2 \, y_2\, (|w_1|^2 - w_2 w_1) = - m_2 \, \beta_1 \,  y_1 \, (|w_2|^2 - w_2 w_1).
\end{equation}

We can now state and prove the following result.

\begin{Theorem} \label{Prop:hyperbolic-equilibria-1}
Consider 2 point particles of masses $m_1, m_2>0$ moving in $\mathbb H_R^2$. Then some necessary and sufficient conditions for the existence of a hyperbolic relative equilibrium as a solution of system \eqref{eq:motionkleingeo} with $n=2$ are that one particle moves along a non-geodesic half line, while the other particle moves along another non-geodesic half line, both half lines converging to the origin of the coordinates system, such that the supporting lines have slopes of opposite signs that satisfy the relationship
\begin{equation}\label{eq:cond-equil-rel-two-bodies-mass-slop}
\frac{m_1}{m_2}=-\frac{\beta_2}{\beta_1} \frac{y_2}{y_1},
\end{equation}
and that at every time instant there is a geodesic half circle centered at the origin of the coordinate system on which both particles are located.
\end{Theorem}
\begin{proof}
If we take the real and imaginary parts of equation (\ref{eq:twopolinomialkleinslopes}),
we obtain the equations
\begin{equation}\label{eq:cond-equil-rel-two-bodies}
m_1 \beta_2 y_2 |w_1|^2 - m_1 \beta_2 y_2 {\rm Re} ( w_2 w_1) = - m_2 \beta_1  y_1  |w_2|^2 + m_2  \beta_1  y_1 {\rm Re} (w_2 w_1),
\end{equation}
\begin{equation}\label{eq:cond-equil-rel-two-bodies2}
-m_1  \beta_2 y_2 {\rm Im}  (w_2 w_1) =  m_2  \beta_1   y_1  {\rm Im}  (w_2 w_1). 
\end{equation}
From equation (\ref{eq:cond-equil-rel-two-bodies2}), it follows that 
$$
-m_1 \beta_2 y_2 = m_2 \beta_1  y_1,
$$
so condition \eqref{eq:cond-equil-rel-two-bodies-mass-slop} must be satisfied. Equation \eqref{eq:cond-equil-rel-two-bodies} implies that
$$
m_1  \beta_2 y_2 |w_1|^2 = - m_2 \beta_1  y_1 |w_2|^2,
$$
which, by condition \eqref{eq:cond-equil-rel-two-bodies-mass-slop}, is equivalent to $|w_1|^2=|w_2|^2$. This fact proves that, at every time instant, there must exist a geodesic half circle centered at the origin of the coordinate system on which the particles are located. This remark completes the proof. 
\end{proof}

\begin{Remark}
Notice that, in the above result, if the masses are equal, $m_1=m_2>0$, then both half lines have slopes equal in absolute value, but of opposite sign. Indeed, 
since, at every time instant, the particles must be located on the same geodesic half circle centered at the origin of the coordinate system, the slopes $\beta_1$ and $\beta_2$ of the half lines along which they move satisfy $\beta_1= - \beta_2$ if and only if $y_1=y_2$. Therefore, from equation (\ref{eq:cond-equil-rel-two-bodies-mass-slop}),
we have that $m_1=m_2$ if and only if $\beta_1= -\beta_2$. 
\end{Remark}

\subsubsection{The case $n = 3$} \label{sec:casehyp_2}

We will next study the case of 3 bodies in the Poincar\'e upper half plane, $\mathbb H_R^2$, with masses $m_1, m_2, m_3>0$. In this context, the system of algebraic
equations (\ref{eq:condrationalsystemklein1a}) becomes
\begin{equation}
\label{eq:threerationalsystemkleinr1}
\frac{R (w_1+\bar{w}_1) \, w_1}{8(w_1- \bar{w}_1)^4}
\end{equation}
$$
=\frac{m_2(w_2-\bar{w}_2)^2(w_1-w_2)(\bar{w}_2-w_1)}{\{[(w_1+\bar{w}_1)(w_2+\bar{w}_2)-2(|w_2|^2+|w_1|^2)]^2-(w_1-\bar{w}_1)^2
(w_2-\bar{w}_2)^2\}^{3/2}} 
$$
$$
+\frac{m_3(w_3-\bar{w}_3)^2(w_1-w_3)(\bar{w}_3-w_1)}{\{[(w_1+\bar{w}_1)(w_3+\bar{w}_3)- 2(|w_3|^2+|w_1|^2)]^2-(w_1-\bar{w}_1)^2(w_3-\bar{w}_3)^2\}^{3/2}},
$$

\begin{equation}
\label{eq:threerationalsystemkleinr2}
\frac{R (w_2+\bar{w}_2) \, w_2}{8(w_2- \bar{w}_2)^4}
\end{equation}
$$
=\frac{m_1(w_1-\bar{w}_1)^2(w_2-w_1)(\bar{w}_1-w_2)}{\{[(w_1+\bar{w}_1)(w_2+\bar{w}_2)-2(|w_2|^2+|w_1|^2)]^2-(w_1-\bar{w}_1)^2
(w_2-\bar{w}_2)^2\}^{3/2}} 
$$
$$
+\frac{m_3(w_3-\bar{w}_3)^2(w_2-w_3)(\bar{w}_3-w_2)}{\{[(w_2+\bar{w}_2)(w_3+\bar{w}_3)- 2(|w_2|^2+|w_3|^2)]^2-(w_2-\bar{w}_2)^2(w_3-\bar{w}_3)^2\}^{3/2}}, 
$$

\begin{equation}
\label{eq:threerationalsystemkleinr3}
\frac{R (w_3+\bar{w}_3) \, w_3}{8(w_3- \bar{w}_3)^4}
\end{equation}
$$
=\frac{m_1(w_1-\bar{w}_1)^2(w_3-w_1)(\bar{w}_1-w_3)}{\{[(w_1+\bar{w}_1)(w_3+\bar{w}_3)-2(|w_1|^2+|w_3|^2)]^2-(w_1-\bar{w}_1)^2
(w_3-\bar{w}_3)^2\}^{3/2}} 
$$
$$
+\frac{m_2(w_2-\bar{w}_2)^2(w_3-w_2)(\bar{w}_2-w_3)}{\{[(w_2+\bar{w}_2)(w_3+\bar{w}_3)- 2(|w_2|^2+|w_3|^2)]^2-(w_2-\bar{w}_2)^2(w_3-\bar{w}_3)^2\}^{3/2}}.
$$

We assume that the particles $m_1$ and $m_3$ move along non-geodesic half lines and that $m_2$ moves along the geodesic vertical half line, such that, at every time instant, there is a geodesic half circle centered at the orgin of the coordinate system on which all particles are located. In other words, if $w=(w_1,w_2,w_3)$ represents the configuration of the system, we have 
$$
w_2+\bar{w}_2=0\ \ {\rm and}\ \ |w_1|=|w_2|=|w_3|.
$$
We can also write that 
$$
w_1(t)=w_1(0)e^{-t/2},\ \  w_2(t)=w_2(0)e^{-t/2}, \ \ w_3(t)=w_3(0)e^{-t/2},
$$
with $|w_1(0)|=|w_3(0)|=1$ and $w_2(0)=i$. The latter conditions are not restrictive,
since the only requirement for the initial conditions is to lie on the half lines on which the particles are assumed to move. Substituting the above forms of $w_1, w_2, w_3$ into
equations \eqref{eq:threerationalsystemkleinr1}, \eqref{eq:threerationalsystemkleinr2}, \eqref{eq:threerationalsystemkleinr3}, the factors $e^{-t/2}$ get cancelled, and after redenoting $w_1(0), w_2(0), w_3(0)$ by $w_1, w_2, w_3$, respectively, we obtain the following equations:
\begin{equation}\label{try1}
\frac{R (w_1+\bar{w}_1) \, w_1}{8(w_1- \bar{w}_1)^4} =
\frac{m_2(w_1^2+1)}{2[4+(w_1-\bar{w}_1)^2]^{3/2}} 
\end{equation}
$$
+\frac{m_3(w_3-\bar{w}_3)^2(w_1-w_3)(\bar{w}_3-w_1)}{\{[(w_1+\bar{w}_1)(w_3+\bar{w}_3)- 4]^2-(w_1-\bar{w}_1)^2(w_3-\bar{w}_3)^2\}^{3/2}},
$$ 

\medskip

\begin{equation}\label{try2}
0=
\frac{m_1(w_1-\bar{w}_1)^2 (w_1 +\bar{w}_1)}{[4+(w_1-\bar{w}_1)^2]^{3/2}}
+ \frac{m_3(w_3-\bar{w}_3)^2(w_3+\bar{w}_3)}{[4+(w_3-\bar{w}_3)^2]^{3/2}}, 
\end{equation}

\medskip

\begin{equation}\label{try3}
\frac{R (w_3+\bar{w}_3) \, w_3}{8(w_3- \bar{w}_3)^4} =
\frac{m_2(w_3^2+1)}{2[4+(w_3-\bar{w}_3)^2]^{3/2}}
\end{equation}
$$
+\frac{m_1(w_1-\bar{w}_1)^2(w_3-w_1)(\bar{w}_1-w_3)}{\{[(w_1+\bar{w}_1)(w_3+\bar{w}_3)-4]^2-(w_1-\bar{w}_1)^2(w_3-\bar{w}_3)^2\}^{3/2}}.
$$

Equation (\ref{try2}) can take place only if
$w_1 +\bar{w}_1$ and $w_3 +\bar{w}_3$ have opposite signs, so we can rewrite this
equation as
\begin{equation}\label{new1}
 \frac{m_1(w_1-\bar{w}_1)^2 |w_1 +\bar{w}_1|}{[2 (w_1+\bar{w}_1)]^{3}} =
\frac{m_3(w_3-\bar{w}_3)^2|w_3+\bar{w}_3|}{[2 (w_3+\bar{w}_3)]^{3}}.
\end{equation}

To express $w_1$ and $w_3$ in terms of the angles the half lines make with
the horizontal axis, we put
$$
w_1=e^{i\theta_1}=\cos\theta_1+i\sin\theta_1\ \ {\rm and}\ \ w_2=e^{i\theta_2}=\cos\theta_2+i\sin\theta_2,
$$
with $\theta_1,\theta_3\in(-\pi/2,0)\cup(0,\pi/2)$. Then equation \eqref{new1} becomes 
\begin{equation}\label{tan}
m_1\tan^2\theta_1=m_3\tan^2\theta_3,
\end{equation}
which shows what relationship exists between the masses and the angles of the
non-geodesic half lines along which the corresponding particles move.

We can now state and prove the main result of this section, which shows that Eulerian relative equilibria for which one body moves along a geodesic exist only if the masses moving on non-geodesic curves equidistant from the geodesic are equal and those curves are on opposite parts of the geodesic and at the same distance from it.

\begin{Theorem} \label{eq:hyper-rel-equi-3-body}
Consider 3 point particles of masses $m_1, m_2, m_3>0$ moving in $\mathbb H_R^2$. Assume that $m_1$ and $m_3$ move along non-geodesic half lines emerging from the origin of the coordinate system at angles $\theta_1$ and $\theta_3$, respectively, and that $m_2$ moves along the geodesic vertical half line. Moreover, at every time instant, there is a geodesic half circle on which all 3 bodies are located, and the motion of the particles is given by the function $w=(w_1,w_2,w_3)$. Then $w$ is a hyperbolic relative equilibrium that is a solution of system
\eqref{eq:motionkleingeo} with $n=3$ if and only if $\theta_1=-\theta_3$ and $m_1=m_3$, with $\theta_1,\theta_3\in(-\pi/2,0)\cup(0,\pi/2)$.
\end{Theorem}
\begin{proof}
We already showed that for $w_1=e^{i\theta_1}, w_2=i,$ and $w_3=e^{i\theta_3}$, equation \eqref{try2} takes the form \eqref{tan}. With the same substitutions, equations 
\eqref{try1} and \eqref{try3} become, respectively, 
\begin{equation}
\label{one}
\frac{R \cos \theta_1}{\sin^4 \theta_1} = \frac{8 m_2}{\cos^2 \theta_1}+
\frac{8 m_3 \sin^2 \theta_3}{(\cos \theta_3-\cos \theta_1)^2},
\end{equation}
\begin{equation}
\label{two}
\frac{R \cos \theta_3}{\sin^4 \theta_3} = \frac{8 m_1 \sin^2 \theta_1}{(\cos \theta_1-\cos \theta_3)^2}+
\frac{8 m_2 }{\cos^2 \theta_3}. 
\end{equation}
If we divide equation \eqref{one} by $\cos^2 \theta_3$ and equation \eqref{two} by $\cos^2 \theta_1$, using relation \eqref{tan} we obtain that
\begin{equation}
\label{sin-cos}
\frac{ \cos^3 \theta_1}{\sin^4 \theta_1} = \frac{\cos^3 \theta_3}{\sin^4 \theta_3}. 
\end{equation}
Consider now the function 
$$
f\colon(-\pi/2,0)\cup(0,\pi/2),\ \ \ f(x)=\frac{\cos^3x}{\sin^4x},
$$
which is obviously even. It is easy to see that $f$ is increasing in the interval
$(-\pi/2,0)$ and decreasing in the interval $(0,\pi/2)$. Therefore equation \eqref{sin-cos}
has solutions if and only if $\theta_1=\pm\theta_3$. Since $\theta_1=\theta_3$ induces
a collision configuration, which is a singularity, the only possible solution is $\theta_1=-\theta_3$. The fact that $m_1=m_3$ follows now from equation \eqref{tan}. This remark completes the proof.
\end{proof}

\subsection{Parabolic relative equilibria}
\label{subsec:parabolic-relative eq}

In this section we will study the  relative equilibria associated to the subgroup
\[ \phi_2 (t)=\exp (t X_2)=
\left(\begin{array}{cc}
    1 & t \\
    0  & 1  \\
    \end{array}\right),
\]
generated by the Killing vector field $X_2$ and which defines the one-parametric family of acting M\"obius transformations
\[ f_{2} (w,t) = w (t)+t ,\]
in the upper half plane $\mathbb{H}^2_R$. These orbits correspond to parabolic
relative equilibria, and we will show that they do not exist in $\mathbb H_R^2$.

Let $\zeta =(\zeta_1,\dots, \zeta_n)$, with $\zeta_k (t) = w_k (t) + t$, be the action orbit for $w=(w_1,\dots, w_n)$, which is a solution  of the equations of motion (\ref{eq:motionkleingeo}). Then
$$
\dot{\zeta}_k = \dot{w}_k+1 \ \ {\rm and}\ \ 
\ddot{\zeta}_k = \ddot{w}_k, \ k=1,\dots, n,
$$
therefore $\zeta$ is  also a solution of system (\ref{eq:motionkleingeo}) if and only if
$$
m_k\ddot{\zeta}_k=\frac{2m_k\dot{\zeta}_k^2}{\zeta_k-\bar{\zeta}_k}-
\frac{(\zeta_k-\bar{\zeta}_k)^2}{2R^2}\frac{\partial V_R}{\partial\bar{\zeta}_k},\ k=1,\dots,n,
$$
which can be written as
$$
 m_k \ddot{w}_k  = \frac{2 m_k (\dot{w}_k +1)^2}{ w_k- \bar{w}_k} -
 \frac{ (w_k- \bar{w}_k)^2}{2 R^2} \, \frac{\partial V_R}{\partial \bar{w}_k} \,\frac{d \bar{w}_k}{d \bar{\zeta}_k}, \ k=1,\dots,n, 
$$
which, since $\displaystyle \frac{d \bar{w}_k}{d \bar{\zeta}_k}=1$, is the same as
\begin{equation}\label{sol-zeta}
 m_k \ddot{w}_k = \frac{2 m_k ( \dot{w}_k^2+ 2 \dot{w}_k +1 ) }{ w_k- \bar{w}_k } - \frac{(w_k-\bar{w}_k )^2}{2 R^2} \,
\frac{\partial V_R}{\partial \bar{w}_k}, \ k=1,\dots,n.
\end{equation}
But since $w$ is a solution of system (\ref{eq:motionkleingeo}), we also have that
\begin{equation}\label{sol-w}
m_k\ddot{w}_k=\frac{2m_k\dot{w}_k^2}{w_k-\bar{w}_k}-
\frac{(w_k-\bar{w}_k)^2}{2R^2}\frac{\partial V_R}{\partial\bar{w}_k},\ k=1,\dots,n,
\end{equation}
Comparing now equations \eqref{sol-zeta} and \eqref{sol-w}, we obtain that
\begin{equation} \label{eq:condiklein2}
2\dot{w}_k = -1,\ k=1,\dots,n,
\end{equation}
which holds if and only if
\begin{equation}\label{eq:principalcondklein2}
   w_k (t) = - \frac{t}{2}+ w_k(0),\ k=1,\dots,n,
\end{equation}
where $w_k(0),\ k=1,\dots, n,$ are initial conditions. Consequently, a necessary condition for the particles $m_1,\dots, m_n$ to form a relative equilibrium associated to the Killing vector field $X_2$ is that they move along horizontal straight lines in $\mathbb{H}^2_R$ passing through $w_k(0), \ k=1,\dots, n$. In terms of the Poincar\'e  disk $\mathbb{D}^2_R$, these orbits correspond to the parametric curves
$$
z_k(t)= \frac{-R [-\frac{t}{2}+ w_k(0)]+ iR^2}{ - \frac{t}{2}+ w_k(0)+ iR},\ k=1,\dots,n,
$$
which start at the point $(-R,0)$, as $t \to -\infty$, and end at the same point $(-R,0)$,
as $t\to\infty$.
These curves have the same topology as the boundary circle of $\mathbb D_R^2$.
In terms of the hyperbolic sphere $\mathbb{L}_R^2$, these lines are the parabolas obtained by intersecting $\mathbb{L}_R^2$ with a plane orthogonal to the rotation axis, the line $y=0, z=x$.

We can now state and prove the following result.

\begin{Theorem}\label{thm:existenceklein2} Consider $n\ge 2$ point particles
of masses $m_1,\dots, m_n>0$ moving in $\mathbb{H}^2_{R}$. Then a necessary and sufficient  condition for  the function $w=(w_1,\dots, w_n)$ to be a solution of system
\eqref{eq:motionkleingeo} that is a relative equilibrium associated to the Killing vector field $X_2$ is that the coordinate functions satisfy the equations
\begin{equation} \label{eq:condrationalsystemklein2b}
-\frac{R}{4(w_k- \bar{w}_k)^4} = \sum_{\substack{j=1\\ j\ne k}}^n
\frac{m_j
(\bar{w}_j-w_j)^2(w_k-w_j)(\bar{w}_j-w_k)}{[\tilde{\Theta}_{3,(k,j)}(w,
\bar{w})]^{3/2}},
\end{equation}
$k=1,\dots, n$, where
\begin{equation}\label{eq:singularklein2c}
\tilde{\Theta}_{3,(k,j)}(w, \bar{w}) = [(\bar{w}_k+w_k)(\bar{w}_j+w_j)-2(|w_k|^2+|w_j|^2)]^2 
\end{equation}
$$
-(\bar{w}_k-w_k)^2(\bar{w}_j-w_j)^2, \ k,j\in\{1,\dots,n\},\ k\ne j.
$$
\end{Theorem}
\begin{proof}
We saw that relative equilibria associated with the Killing vector field $X_2$ must
satify equations (\ref{eq:principalcondklein2}), which imply that
\begin{equation} \label{harmonicklein2}
\ddot{w}_k = 0,\ k=1,\dots,n.
\end{equation}
Therefore, from equations \eqref{sol-w}, we can conclude that the coordinates of a  relative equilibrium $w$ satisfy the equations
\begin{equation} \label{eq:condrationalsystemklein2}
\frac{m_k R^2}{(w_k- \bar{w}_k)^3} = \frac{\partial V_R}{\partial \bar{w}_k},
\ \ k=1,\dots, n.
\end{equation}
If we now compare the above equations to the expressions \eqref{eq:gradmet2} of
$\frac{\partial V_R}{\partial \bar{w}_k}, \ k=1,\dots, n$, we obtain the desired
relationships \eqref{eq:singularklein2c}. This remark completes the proof.
\end{proof}

\begin{Definition} \label{def:paraprelequil} We will call {\it parabolic relative equilibria} the solutions of system \eqref{eq:motionkleingeo} in $\mathbb{H}^2_{R}$ that satisfy
equations \eqref{eq:condrationalsystemklein2b}.
\end{Definition}

Notice that equations (\ref{eq:condiklein2}) provide the velocities of the particles in case they form a parabolic relative equilibrium. However, as we will further prove, parabolic relative equilibria do not exist in the curved $n$-body problem. The following statement generalizes a result obtained in \cite{Diac} for curvature $\kappa=-1$. Using the same idea as in \cite{Diac}, this result was generalized in \cite{Diacu3} to the 3-dimensional case.

\begin{Theorem}\label{theo:non-parabolic-relative-equilibria}
In the curved $n$-body problem with negative curvature there are no parabolic relative equilibria.
\end{Theorem}
\begin{proof}
Using the notation
$$
w_k=a_k+ib_k \ \ {\rm and} \ \  w_j=a_j+ib_j,
$$
the real and imaginary part of equations \eqref{eq:condrationalsystemklein2b} become
\begin{equation} \label{eq:condrationalsystemklein2d}
-\frac{R}{32 b_k^4} = \sum_{\substack{j=1\\ j\ne k}}^n
\frac{2  m_j b_j^2 [(a_k-a_j)^2+(b_j^2-b_k^2)]}{[\tilde{\Theta}_{3,(k,j)}(w_0,
\bar{w}_0)]^{3/2}}, \ k=1,\dots, n,
\end{equation}
\begin{equation}\label{eq:condrationalsystemklein2d1}
0 = \sum_{\substack{j=1\\ j\ne k}}^n
\frac{4 m_j  b_j^2  b_k  (a_k-a_j)}{[\tilde{\Theta}_{3,(k,j)}(w_0,
\bar{w}_0)]^{3/2}}, \ k=1,\dots, n.
\end{equation}

Since $b_k>0, \ k=1,\dots, n$, equations \eqref{eq:condrationalsystemklein2d1} hold for any $k, j\in\{1,\dots,n\}$, with $k\ne j$, if and only if $a_{k}=a_{j}$. This fact implies that all the particles are located on the same vertical line. Without loss of generality, we can assume that they are on the vertical half line $x=0,\ y>0$. Therefore  $w_{k}=b_{k}i$ and $w_{j}=b_{j} i$.
When we substitute these values into equation
(\ref{eq:condrationalsystemklein2b}) we obtain
\begin{equation} \label{eq:condrationalsystemklein2c}
-\frac{R}{32 b_{k}^4} = \sum_{\substack{j=1\\ j \neq k}}^n \frac{m_jb_{j}^2(b_j^2-b_k^2)}{|b_{k}^2-b_{j}^2|^3}, \ k=1,\dots,n.
\end{equation}
Since the particles do not collide, we can assume, without loss of generality, that
$$0<b_1<\dots<b_n.$$ Then, for $k=1$, we can conclude from \eqref{eq:condrationalsystemklein2c} that
$$
-\frac{R}{32 b_{1}^4} = \sum_{j=2}^n \frac{m_jb_{j}^2(b_j^2-b_1^2)}{|b_{1}^2-b_{j}^2|^3}.
$$
But the left hand side of this equation is negative, whereas the right hand side is positive. This contradiction completes the proof.
\end{proof}

\subsection*{Acknowledgments}
Florin Diacu acknowledges the partial support of an NSERC Discovery Grant, whereas Ernesto P\'erez-Chavela and J.\ Guadalupe Reyes Victoria acknowledge the partial support received from Grant 128790 provided by
CONACYT of M\'exico.

\end{document}